\documentclass[a4paper,11pt]{amsart}
\usepackage{tikz}
\usepackage{amsfonts,amsmath,amssymb}
\usepackage{color}
\usepackage{algorithmic}
\usepackage{algorithm}
\usepackage{comment}
\usepackage{graphics}
\usepackage[hidelinks]{hyperref}

\usepackage{listings}
\usepackage{mathtools}
\usepackage{esint}

\theoremstyle{definition}

\theoremstyle{lemma}

\theoremstyle{remark}

\numberwithin{theorem}{section}
\numberwithin{equation}{section}
\numberwithin{table}{section}
\numberwithin{figure}{section}

\usepackage{caption}
\usepackage{subcaption}
\usepackage{graphics}

\definecolor{aliceblue}{rgb}{0.94, 0.97, 1.0}

\newcommand{\ddiv}{\operatorname{div}}
\newcommand{\calT}{\mathcal{T}}
\newcommand{\R}{\mathbb{R}}
\newcommand{\N}{\mathbb{N}}

\newcommand{\opnorm}[1]{{\left\vert\kern-0.25ex\left\vert\kern-0.25ex\left\vert #1 \right\vert\kern-0.25ex\right\vert\kern-0.25ex\righ\vert}}

\newcommand\dx{\,\text{d}x}

\newcommand{\opL}{\mathfrak{L}}

\newcommand{\frakA}{\mathfrak{A}}

\newcommand{\tVh}{\tilde{V}_h}
\newcommand{\Vh}{V_h}
\newcommand{\V}{H^1_0(D)}

\newcommand{\LL}{L^2(D)}

\newcommand\calP{\mathcal{P}}
\newcommand\calQ{\mathcal{Q}}

\newcommand\calI{\mathcal{I}}
\newcommand\calD{\mathcal{D}}
\newcommand{\eps}{\varepsilon}

\newcommand{\W}{\mathcal{W}}

\newcommand{\Nb}{\mathtt{N}}

\definecolor{darkgreen}{rgb}{0.09, 0.45, 0.27}
\definecolor{myOrange}{rgb}{0.85000,0.32500,0.09800}

\newcommand\corr[1]{{#1}}

\begin{document}
%
%
\title[Operator Compression with Deep Neural Networks]{
Operator Compression with Deep Neural Networks}
\author[F.~Kr\"opfl, R.~Maier, D.~Peterseim]{Fabian Kr\"opfl$^\dagger$, Roland Maier$^\ddagger$, Daniel Peterseim$^\dagger$$^\star$}
\address{${}^{\dagger}$ Institute of Mathematics, University of Augsburg, Universit\"atsstr.~12a, 86159 Augsburg, Germany}
\email{fabian.kroepfl@uni-a.de \\ daniel.peterseim@uni-a.de}
\address{${}^{\star}$ Centre for Advanced Analytics and Predictive Sciences (CAAPS), University of Augsburg, Universit\"atsstr.~12a, 86159 Augsburg, Germany}
\address{${}^{\ddagger}$ Institute of Mathematics, Friedrich Schiller University Jena, Ernst-Abbe-Platz 2, 07743 Jena, Germany}
\email{roland.maier@uni-jena.de}
\thanks{The work of F.~Kr\"opfl and D.~Peterseim is part of a project that has received funding from the European Research Council (ERC) under the European Union's Horizon 2020 research and innovation programme (Grant agreement No.~865751 - \emph{Computational Random Multiscale Problems}). R.~Maier acknowledges support by the G\"oran Gustafsson Foundation for Research in Natural Sciences and Medicine.}
\date{\today}

\maketitle
%
%
\begin{abstract}
This paper studies the compression of partial differential operators using neural networks. We consider a family of operators, parameterized by a potentially high-dimensional space of coefficients that may vary on a large range of scales. Based on existing methods that compress such a multiscale operator to a finite-dimensional sparse surrogate model on a given target scale, we propose to directly approximate the coefficient-to-surrogate map with a neural network. We emulate local assembly structures of the surrogates  and thus only require a  moderately sized network that can be trained efficiently in an offline phase. This enables large compression ratios and the online computation of a surrogate based on simple forward passes through the network is substantially accelerated compared to classical numerical upscaling approaches. We apply the abstract framework to a family of prototypical second-order elliptic heterogeneous diffusion operators as a demonstrating example.
\end{abstract}


%
%
\section{Introduction}\label{sec: intro}
The remarkable success of machine learning technology, especially deep learning, in classical AI disciplines such as image recognition and natural language processing has led to an increased research interest in leveraging the power of these approaches in other science and engineering disciplines over the last years.
In the field of numerical modeling and simulation, promising approaches are emerging that try to integrate machine learning algorithms and traditional physics-based approaches, combining the advantages of the data-driven regime with known physics and domain knowledge. In this spirit, many different approaches to approximating solutions of partial differential equations (PDEs) with neural networks have been proposed, for example so-called Physics-Informed Neural Networks (PINNs)~\cite{RaiPK19}, the Deep Galerkin method~\cite{SirS18} or the Deep Ritz Method~\cite{WeiY18}.
It has become evident that the strategy of using neural networks as ansatz functions for the approximation of a PDE's solution is especially advantageous for high-dimensional problems that are outside the reach of classical mesh-based methods~\cite{HanJ20,HanJW18,WeiHJ17}. For some classes of PDEs, e.g., Kolmogorov PDEs and semilinear heat equations, it has even been proven that neural networks break the curse of dimensionality~\cite{BerDG20,HutJKN20}.

In this spirit, we strongly believe that the strength of neural networks lies in scenarios where one deals with a whole family of PDEs rather than one single equation, for example in the context of so-called parametric PDEs, i.e., settings where a family of partial differential operators parameterized by some coefficient is considered, see, e.g., \cite{SchZ19}. This  is particularly true for multiscale problems, where one is interested in computing coarse-scale surrogates for problems involving a range of scales that cannot be resolved in a direct numerical simulation.

In this paper, we study the problem of approximating a coefficient-to-surrogate map with a neural network in a very general setting of parameterized PDEs with arbitrarily rough coefficients that may vary on a microscopic scale. In other words, we are not trying to directly approximate the parameter-to-solution map, but rather compress the fine-scale information contained in the continuous operator to a finite-dimensional sparse object that is able to replicate the effective behavior of the solution on a macroscopic scale of interest even in the presence of unresolved oscillations of the underlying coefficient.

The output surrogate models are based on the idea of modern numerical homogenization techniques such as Localized Orthogonal Decomposition \cite{MalP14,Pet16,MP20}, Gamblets \cite{Owh17}, Rough Polyharmonic Splines \cite{OwhZB13}, the Multiscale Finite Element Method~\cite{Hou1997,EfeH09}, or the Generalized Finite Element Method \cite{BabL11,EfeGW11}; see \cite{AltHP21} and the references therein for a comprehensive overview. These methods have demonstrated high performance in many relevant applications such as porous media flow or wave scattering in heterogeneous media to mention only a few. In particular, they typically do not require explicit assumptions on the existence of lower-dimensional structures in the underlying family of PDE coefficients and yield sparse system matrices that ensure uniform approximation properties of the resulting surrogate. Moreover, the computation of the system matrices mimics the standard assembly procedure from finite element theory, consisting of the generation of local system matrices and their combination by local-to-global mappings, which is exploited to reduce the size of the network architecture and its complexity considerably.

The possibility of fast computation of the surrogates has high potential for multi-query problems, such as in uncertainty quantification, and time-dependent or inverse multiscale problems, which require the computation of surrogates for many different a priori unknown coefficients. Though the aforementioned numerical homogenization methods lead to accurate surrogates for the whole class of coefficients, their computation requires the resolution of all scales locally which marks a severe limitation when it has to be performed many times for the solution of a multi-query problem.  
There have been attempts to tackle this problem, but the results so-far are only applicable to small perturbation regimes~\cite{HelKM20, MV21} or settings where the parameterization fulfills additional smoothness requirements~\cite{AbdH15}. 

To overcome this problem, we propose to learn the whole nonlinear coefficient-to-surrogate map from a training set consisting of pairs of coefficients and their corresponding surrogates with a deep neural network. In other words, we are combining the domain knowledge from numerical homogenization with a data-driven deep learning approach by essentially learning a numerical homogenization method from data. To this end, we propose using an offline-online approach. In the offline phase, the neural network is trained based on data generated with existing numerical homogenization techniques. In the online phase, the compression of previously unseen operators can then be reduced to a simple forward pass of the neural network, which eliminates  the computational bottleneck encountered in multi-query settings.

Our method is conceptually different from existing approaches that try to integrate ideas from homogenization with neural networks. In~\cite{ArbBSRK20} for example, the authors propose to learn a homogenized PDE from simulation data by linking deep learning with an equation-free multiscale approach. Other papers in the context of uncertainty quantification suggest training a neural network to identify suitable multiscale basis functions for the finite volume method given a porous random medium~\cite{ChaE18,PadZ21}. In~\cite{GhaS19}, the authors consider the problem of elasticity with history-dependent material properties, {where a} recurrent deep neural network connects microscopic and macroscopic material parameters. In Deep Multiscale Model Learning~\cite{WanCCEW20}, {learning techniques are used to predict the evolution from one time step to another within a given coarse multiscale space. The goal of this approach is to obtain a reasonable coarse operator for the successive approximation of a time-dependent PDE.} 
Furthermore, several experimental and theoretical works on the approximation of the coefficient-to-solution map \cite{GeiPRSK20,BhaHKS20,GaoSW21,KutPRS19} or other quantities of interest such as the ground state energy in Schr\"odinger equations \cite{KhoLY17} by deep neural networks have been published in the context of parametric PDEs.

This paper is structured as follows: in Section \ref{sec:framework}, we introduce and motivate the abstract framework for a very general class of linear differential operators. After that, we study the problem of elliptic homogenization as an example of how to apply the general methodology in practice. In Section \ref{sec:numexp}, we conduct numerical experiments that show the feasibility of our ideas developed in the previous two sections. We conclude this work with an outlook on further research questions.

%
%
\section{Abstract Framework}\label{sec:framework}

In this section, we describe the general abstract problem of finding discrete compressed surrogates to a family of differential operators that allow us to satisfactorily approximate the original operators on a target scale of interest, given only the underlying coefficients but \emph{not} a high resolution representation of the operators. We elaborate on how to speed up the online computation of those compressed representatives using deep neural networks after an initial offline training phase.

\subsection{Setting}\label{subsec:setting}
Let $D \subseteq \R^d$, $d\in\{1,2,3\}$ be a bounded Lipschitz domain and $H^1_0(D)$ the Sobolev space of $L^2$-functions with weak first derivatives in $L^2(D)$ that vanish on the boundary of $D$. We write $H^{-1}(D)$ for the dual space of $H^1_0(D)$ and $\langle\cdot,\cdot\rangle$ for the duality pairing between $H^{-1}(D)$ and $H^1_0(D)$.
Consider a family of linear differential operators
\begin{equation*}
\boldsymbol{\opL} := \{\opL_A \colon H^1_0(D) \to H^{-1}(D)\ \vert\  A\in\mathfrak{A}\},
\end{equation*}
that is parameterized by some class $\mathfrak{A}\subseteq L^\infty(D)$ of admissible coefficients. We emphasize, that we do not pose any assumptions on the structure of the coefficients $A\in\mathfrak{A}$ such as periodicity or scale separation and explicitly allow for arbitrarily rough coefficients that may vary on a continuum of scales up to some microscale $\varepsilon \ll \text{diam}(D)$. We assume that for every $A\in\mathfrak{A}$, the associated operator $\opL_A$ is symmetric ($\langle \opL_Au,v\rangle = \langle u,\opL_A v \rangle)$, local ($\langle \opL_A u,v \rangle = 0$ if $u$ and $v$ have disjoint supports) and bijective.

Bijectivity implies that for any given $A\in\mathfrak{A}$ and $f\in H^{-1}(D),$ there exists a unique $u\in H^1_0(D)$ that solves the equation
\begin{equation}\label{eq:abstractmodel}
\opL_A u = f
\end{equation} 
in a weak sense, i.e., the solution satisfies
\begin{equation}\label{eq:abstractmodelweak}
\langle \opL_A u, v \rangle = \langle f, v \rangle \quad \text{for all } v\in\V.
\end{equation} 
For given problem data $A$ and $f$, we are interested in computing an approximation to $u$ on some target scale in reasonable time. 

\subsection{Discretization}\label{subsec:discretization}
In order to be able to solve this problem computationally, we choose a finite-dimensional subspace $V_h\subseteq\V$ of dimension $m = \mathrm{dim}(\Vh)$. As a standard example, one could take $\Vh$ to be a classical finite element space based on some mesh $\calT_h$ with characteristic mesh size $h$ and approximate \eqref{eq:abstractmodelweak} with a Galerkin method. However, in the very general setting with $A$ possibly having fine oscillations on a scale that is not resolved by the mesh size $h$, this approach leads to unreliable approximations of $u$. Then again, the resolution of these fine-scale features can be prohibitively expensive in terms of computational resources if $\varepsilon$ is very small. 
Note that resolution here may mean that the actual mesh size is significantly smaller than $\varepsilon$, depending on the oscillations of the coefficient and its regularity \cite{BabO00,PetS12}. 
This means that more advanced discretization techniques are required to still obtain reasonable approximations in the unresolved setting. In practice, the challenge is therefore to compress the fine-scale information that is contained in the operator $\opL_A$ to a suitable surrogate $\mathfrak{S}_A$ on the target scale $h$, i.e., the surrogate $\mathfrak{S}_A$ must be chosen in such a way that it is still able to capture the characteristic behavior of the operator $\opL_A$ on the scale of interest. Moreover, we require $\mathfrak{S}_A$ to be a bijection that maps the space $\Vh$ to itself. This ensures that for any $A\in\mathfrak{A}$ and $f\in H^{-1}(D)$ we can find a unique $u_h\in\Vh$ that weakly solves the discretized equation
\begin{equation}\label{eq:abstractmodeldisc}
\mathfrak{S}_A u_h = f_h,
\end{equation}
with $f_h = \mathfrak{M}f$, where $\mathfrak{M}$ is a quadrature-type operator that maps a function in $H^{-1}$ to an appropriate approximation in $\Vh$. Problem \eqref{eq:abstractmodeldisc} needs to be understood as finding $u_h$ that satisfies
\begin{equation*}
\langle \mathfrak{S}_A u_h, v_h \rangle = \langle f_h, v_h \rangle \quad \text{for all } v_h\in \Vh.
\end{equation*} 
The choice of the surrogate is obviously highly dependent on the problem at hand, see for example Section~\ref{subsec:discelliptic} for possible choices in the case of second-order elliptic diffusion operators.

\subsection{Characterization in terms of a system matrix}\label{subsec:matrix}
 We restrict our discussion to choices of surrogates that can be represented by an $m\times m$ system matrix $S_A$ that is often called the \emph{effective} system matrix in the following. We assume that $S_{A}\in\R^{m \times m}$ is of the form $(S_A)_{ij} = \langle \mathfrak{S}_A\lambda_j, \lambda_i \rangle$ for a basis $\lambda_1,\dots,\lambda_m$ of $\Vh$. Note that the basis should be chosen as localized as possible in order for the resulting system matrix to be sparse. The process of operator compression can then be formalized by a compression operator
\begin{equation*}
\mathfrak{C}\colon \mathfrak{A} \to \mathbb{R}^{m\times m} 
\end{equation*}
that maps a given coefficient $A$ to the system matrix $S_{A}$ representing the compressed surrogate $\mathfrak{S}_A$ of the operator $\opL_A$. Once $\mathfrak{C}$ has been evaluated for given $A\in\mathfrak{A}$, the solution to \eqref{eq:abstractmodelweak} can then be approximated with a function $u_h\in \Vh$ for any right-hand side $f\in H^{-1}(D)$ by solving the linear system $S_{A} U = F,$ where $F\in\mathbb{R}^m$ is the vector with entries $F_i : = \langle {\mathfrak{M}}f, \lambda_i \rangle$ and $U\in\mathbb{R}^m$ contains the coefficients of the basis representation $u_h = \sum_{i=1}^m U_i \lambda_i$. 

\subsection{Multi-query scenarios}\label{subsec:multiquery}
For many classes of coefficients $\mathfrak{A}$ and based on the choice of the surrogate, evaluating $\mathfrak{C}$ requires solving local auxiliary problems, during which the finest scale $\varepsilon$ has to be resolved at some point. While this is acceptable if one wants to compress only a few operators in an offline computation, it becomes a major problem once $\mathfrak{C}$ has to be evaluated for many different coefficients $A$ in an online phase, as for example in certain inverse problems, uncertainty quantification, or the simulation of evolution equations with time-dependent coefficients. This motivates a data-driven offline-online approach, where the offline phase consists of training a neural network to approximate the compression operator $\mathfrak{C}$, such that in the subsequent online phase, the evaluation of $\mathfrak{C}$ can be replaced with a simple forward pass through the network, thus eliminating the computational bottleneck. 

\subsection{System matrix decomposition}\label{subsec:matrixdecomp} 
In principle, one could try to {directly} approximate the global operator $\mathfrak{C}$ with a neural network. If the coefficient involves oscillations on some fine scale $\varepsilon$, this would lead to a network architecture with an input layer of size $\mathcal{O}(\varepsilon^{-d})$, an output layer of size $\mathcal{O}(m)$, and possible hidden layers. Particularly for small $\varepsilon$, this leads to very large networks and, thus, requires a huge amount of free parameters and therefore extraordinary amounts of training data and storage space in order to preserve good generalization capabilities.

To reduce the necessary size of the network, one can exploit available information on the compression operator $\mathfrak{C}$ by means of a certain structure in the resulting effective matrices $S_A$. To this end, we think of $S_A$ as a matrix composed of multiple \corr{inflated} sub-matrices, i.e.,
\begin{equation}\label{eq:decomp1}
S_{A} = \sum_{j\in J} \Phi_{j} (S_{A,j}),
\end{equation}
where $J$ denotes some given index set and $S_{A,j}\in\R^{s\times t},\ s,t \ll m,$ are \corr{(typically dense)} local matrices of equal size. The functions
\begin{equation}\label{eq:loc2glob}
\Phi_{j}\colon \R^{s\times t} \rightarrow \R^{m\times m}
\end{equation}
represent local-to-global mappings {inspired by classical finite element assembly processes} as further explained below.
More precisely, there exist index transformations $\pi_j$ and $\varphi_j$, such that 
\begin{equation*}
\Phi_j(S)[\pi_j(k),\varphi_j(l)] = S[k,l],\qquad 1 \leq k \leq s,\, 1 \leq l \leq t,
\end{equation*}
where $M[i_r,i_c]$ denotes the entry of a matrix $M$ in the $i_r^{\hspace*{1pt}\mathrm{th}}$ row and the $i_c^{\hspace*{1pt}\mathrm{th}}$ column. Note that $\pi_j$ and $\varphi_j$ can also map to zero, indicating that the corresponding entry should be disregarded. 
Let us also emphasize that the mappings $\Phi_j$ (and, in turn, the index transformations $\pi_j$ and $\varphi_j$) are completely independent of coefficients $A$ and solely depend on the domain $D$ as well as the geometry of an allotted discretization. \corr{The precise definitions of the maps $\pi_j$ and $\varphi_j$ as well as the index set $J$ are usually determined in a canonical way by the choice of the computational mesh and the compression operator $\mathfrak{C}$. 
In Section~ \ref{subsec:surrogate} below, we present an example of how such mappings may look like.}

Depending on the compression operator $\mathfrak{C}$ and the decomposition \eqref{eq:decomp1}, we can expect that all the local matrices $S_{A,j}$ are created in a similar fashion and only depend on a local sub-sample of the coefficient. This can be understood as a generalization of the assembly process that underlies classical finite element system matrices: these matrices are composed of local system matrices that are computed on each element separately and only require knowledge about the coefficient on the respective element. In that context, the local sub-matrices all have a similar structure and the mapping by the functions $\Phi_j$ leads to overlapping contributions on the global level. 
Going back to the abstract setting, we generalize these properties and assume the existence of a lower-dimensional reduced compression operator
\begin{equation}\label{eq:opcompressionred}
\mathfrak{C}_{\mathrm{red}}\colon \R^r \to \R^{s \times t},
\end{equation} 
such that the contributions $S_{A,j}$ are of the form 
\begin{equation}\label{eq:decomp2}
S_{A,j} = (\mathfrak{C}_{\mathrm{red}} \circ R_j)(A),
\end{equation}
where the operators
\begin{equation}\label{eq:reduction}
R_{j}\colon \mathfrak{A} \rightarrow \R^{r}
\end{equation}
extract $r$ relevant features of a given global coefficient.
In the context of, e.g., finite element matrices, the operators $R_j$ correspond to the restriction of a coefficient to an element-based piecewise constant approximation and $\mathfrak{C}_\mathrm{red}$ incorporates the computation of a local system matrix based on such a sub-sample of the coefficient. 
To achieve a uniform length $r$ of the output for the operators $R_j$, these operators may include artificially introduced zeros depending on the respective geometric configurations (e.g., at the boundary). An example for a quadrilateral mesh in two dimensions is shown in Figure~\ref{fig:boundarypatch}.

The problem of evaluating $\mathfrak{C}$ can now be decomposed into multiple evaluations of the reduced operator $\mathfrak{C}_\mathrm{red}$ that takes the local information $R_j(A)$ of $A$ and outputs a corresponding local matrix as described in~\eqref{eq:decomp2}.
In our setting of a coefficient $A$ that is potentially unresolved by the target scale $h$, evaluating $\mathfrak{C}_{\text{red}}$ is nontrivial and might become a bottleneck in multi-query scenarios as already indicated in Section~\ref{subsec:multiquery} above. In such cases, we propose to approximate the operator $\mathfrak{C}_{\text{red}}$ with a deep neural network 
\begin{equation*}
\Psi(\cdot,\theta): \R^r \rightarrow \mathbb{R}^{s \times t},
\end{equation*}
where $\theta\in\mathbb{R}^p$ is a set of $p$ trainable parameters of moderate size, such that for given $A\in\mathfrak{A}$, the effective system matrix $\mathfrak{C}(A) = S_{A}$ can be efficiently approximated by
\begin{equation*}
\widehat{S}_{A} := \sum_{j\in J} \Phi_j (\Psi(R_j(A),\theta)),
\end{equation*}
which requires just a single forward pass of the minibatch $(R_j(A))_{j\in J}$ through the network. \corr{Note that the approximation $\widehat{S}_A$ possesses the same sparsity structure as the matrix $S_A,$ since the neural network yields only approximations to the local sub-matrices $S_{A,j},$ whereas the assembling process which determines the sparsity structure of the global matrix is determined by the mappings $\Phi_j,$ which are independent of the network $\Psi$.}

We emphasize that a decomposition of $S_A$ as described in \eqref{eq:decomp1}--\eqref{eq:reduction} does not necessarily require a uniform operator $\mathfrak{C}_\mathrm{red}$. If multiple reduced operators are required for such a decomposition, the idea of approximating them by one single neural network can still be applied. It is, however, necessary for the ability of the network to generalize well beyond data seen during training that the reduced operators at least involve certain similarities.

\subsection{Network training}\label{subsec:training} 
In practice, the neural network $\Psi(\cdot,\theta)$ has to be trained in an offline phase from a set of training examples before it can be used for approximating the mapping $\mathfrak{C}_{\mathrm{red}}$.
We propose to draw $N$ global coefficients $(A^{(i)})_{i=1}^N$ from $\mathfrak{A},$ extracting the relevant information $(A_j^{(i)}):=(R_j(A^{(i)}))$ from them and compressing it into the corresponding effective matrices $(S_{A,j}^{(i)})$ with $\mathfrak{C}_{\text{red}}$. This results in a total of  $\vert J\vert \cdot N$ training samples available for the neural network to train on, namely $(A^{(i)}_j,S_{A,j}^{(i)}),\ i = 1,\dots N,\ j \in J.$ In order to learn the parameters of the network, we then minimize the loss functional
\begin{equation}\label{eq:loss}
\mathcal{J}(\theta) = \frac{1}{N\cdot \vert J\vert} \sum_{i=1}^N \sum_{j\in J} \ \frac12 \frac{\| \Psi(A_j^{(i)},\theta) - S_{A,j}^{(i)} \|^2_{\mathbb{R}^{s\times t}}}{\| S_{A,j}^{(i)} \|_{\mathbb{R}^{s\times t}}^2}
\end{equation}
over the parameter space $\mathbb{R}^p$ using iterative gradient-based optimization on minibatches of the training data. This can be very efficiently implemented within modern deep learning frameworks such as \emph{TensorFlow}~\cite{AbaAB+15}, \emph{PyTorch}~\cite{PasG+19}, or \emph{Flux}~\cite{Inn18}, which allow for the automatic differentiation of the loss functional with respect to the network parameters. 

\subsection{Full algorithm}\label{subsec:fullalgo} 
After having established all the conceptual pieces, we now put them together and return to the abstract variational problem \eqref{eq:abstractmodelweak} from the beginning of the section. Suppose that we want to solve \eqref{eq:abstractmodelweak} for a large number of given coefficients $A^{(i)},\ i = 1, \dots, M,$ and a given right hand side $f\in H^{-1}(D)$. For ease of notation, we restrict ourselves to a single right hand side, {which is, however, not necessarily required for our approach}. The proposed procedure is summarized in Algorithm \ref{algo:opcompress} below, divided into the offline and online stage of the method.\\
%
\begin{algorithm}
	\caption{Operator Compression with Neural Network}
	\label{algo:opcompress}
	\begin{itemize}
	 \item[(i)] \textbf{Offline Phase} 
	    \begin{algorithmic}
	 
		\REQUIRE $\text{number of samples }N,\ \text{index set }J, \ \text{operators } R_j, \ \text{error functional } \mathcal{J}$ \\ \hspace*{0.75cm}$\ \ \text{reduced compression operator } \mathfrak{C}_{\text{red}},\ \text{initial neural network } \Psi (\cdot, \theta), $ 
		\ENSURE trained neural network $\Psi (\cdot, \theta)$ \\[1ex]
		\STATE Draw $N$ samples $(A^{(i)})_{i=1}^N$ from $\mathfrak{A}$
		\FOR{$i = 1, \dots, N$} 
		\FOR{$j \in J$}
		\STATE {Extract relevant information from coefficient: $A_j^{(i)} = R_j(A^{(i)})$} 
		\STATE {Compress to local matrix: $S_{A,j}^{\hspace*{1pt} (i)} = \mathfrak{C}_{\text{red}}(A^{(i)}_j)$} 
        \ENDFOR
		\ENDFOR
        \STATE {Train neural network $\Psi(\cdot, \theta)$: update parameters $\theta$ based on data $(A^{(i)}_j, S_{A,j}^{\hspace*{1pt} (i)})_{i,j},$ \\ $i=1,\dots,N,\ j\in J,$ by minimizing the error functional $\mathcal{J}$ over the parameter space $\mathbb{R}^p$ }
 		\RETURN $\Psi(\cdot, \theta)$
		\end{algorithmic}
	\item[(ii)] \textbf{Online Phase}
	\begin{algorithmic}
	\REQUIRE \text{coefficients }$(A^{(i)})_{i=1}^M,\ \text{index set }J,\ \text{operators } R_j,\ \text{trained network }\Psi(\cdot,\theta),$ \\ \hspace*{0.75cm}$\ \ \text{right-hand side vector }F,\ \text{local-to-global mappings } \Phi_{j}$
		\ENSURE coefficient vectors $(U^{(i)})_{i=1}^M$ \\[1ex]
		\FOR{$i = 1, \dots, M$} 
		\FOR{$j \in J$}
		\STATE {Extract relevant information from coefficient: $A_j^{(i)} = R_j(A^{(i)})$} 
        \ENDFOR
		\STATE {Fast compression: $\widehat{S}_{A}^{\hspace*{1pt} (i)} = \sum_{j\in J} \Phi_{j} (\Psi(A_j^{(i)},\theta))$}
			\STATE {Solve: $U^{(i)} = (\widehat{S}_{A}^{\hspace*{1pt} (i)})^{-1} F$}
		\ENDFOR 	
		\RETURN $(U^{(i)})_{i=1}^M$
	\end{algorithmic}
	\end{itemize}	
\end{algorithm}

%
%
\section{Application to Elliptic Homogenization}\label{sec:numhom}
In this section, we specifically consider a family of prototypical elliptic diffusion operators as a demonstrating example of how to apply the abstract framework laid down in Section~\ref{sec:framework} in practice.

\subsection{Setting}\label{subsec:settingelliptic}
From now on let the domain $D$ be polyhedral. We consider the family of linear second-order diffusion operators
\begin{equation*}
\boldsymbol{\opL} := \{-\ddiv(A\nabla \cdot): H^1_0(D) \rightarrow H^{-1}(D)\ \vert\  A\in\mathfrak{A}\},
\end{equation*}
parameterized by the following set of admissible coefficients which possibly encode microstructures, 
\begin{equation}\label{eq:frakAelliptic}
\mathfrak{A} := \left\{ A\in L^\infty(D)\ \vert\ \exists\ 0<\alpha\leq\beta<\infty: \alpha \leq A(x) \leq \beta\ \text{for almost all } x\in D\right\}.
\end{equation}
For the sake of simplicity, we restrict ourselves to scalar coefficients here. Note, however, that also the consideration of matrix-valued coefficients {is not an issue} from a numerical homogenization viewpoint.
We remark that the family of operators $\boldsymbol{\opL}$ fulfills the assumptions of locality and symmetry from the abstract framework. In this setting, the abstract problem \eqref{eq:abstractmodel} amounts to solving the following linear elliptic PDE with homogeneous Dirichlet boundary condition,
\begin{equation*}
\left\{
\begin{aligned}
-\ddiv(A\nabla u) &=f \quad \mathrm{in}\ D,\ \\
u &= 0\quad \mathrm{on}\ \partial D,
\end{aligned} \
\right.
\end{equation*}
which possesses a unique weak solution $u\in\V$ for every $f\in H^{-1}(D)$ {and $A\in\frakA$}. The corresponding counterpart to the weak formulation \eqref{eq:abstractmodelweak}  can be written as: find $u\in\V$ such that  
\begin{equation}\label{eq:modelelliptic}
a_A(u,v) := \int_D A \nabla u \cdot\nabla v\dx = \langle f,v \rangle \quad\text{for all } v \in \V
\end{equation}
by using integration by parts on the divergence term.

\subsection{Discretization and compression}\label{subsec:discelliptic}
Let now $\calT_h$ be a Cartesian mesh with characteristic mesh size $h$ and denote with $Q^1(\calT_h)$ the corresponding space of piecewise bilinear functions. We consider the conforming finite element space $\Vh:=Q^1(\calT_h) \cap \V$ of dimension $m := \text{dim}(\Vh)$. 
Generally, also other types of meshes and finite element spaces could be employed but we restrict ourselves to the above choice for the moment. As we have already mentioned in Section~\ref{subsec:discretization}, if the mesh $\calT_h$ does not resolve the fine-scale oscillations of $A$, approximating $u$ with a pure finite element ansatz of seeking $u_h \in \Vh$ such that
\begin{equation*}
a_A(u_h,v_h) = \langle f,v_h\rangle \quad\text{for all }v_h \in \Vh
\end{equation*} 
will not yield satisfactory results. In a setting where resolving $A$ with the mesh is computationally too demanding, we are therefore interested in suitable choices for a compression operator $\mathfrak{C}$. In particular, we want $\mathfrak{C}$ to produce effective system matrices on the target scale $h$ that can be used to obtain appropriate approximations on this scale. In the following, we briefly comment on possible choices for this operator that are based on the finite element space $\Vh$. 

\subsubsection{Compression by analytical homogenization}\label{sss:anahom}
The idea of analytical homogenization is to replace an oscillating~$A$ by an appropriate homogenized coefficient $A_\mathrm{hom} \in L^\infty(D,\R^{d\times d})$. 
The mathematical theory of homogenization can treat very general non-periodic coefficients in the framework of $G$- or $H$-convergence \cite{Mur78,Spagnolo:1968,Gio75}. However, apart from being non-constructive in many cases, homogenization in the classical analytical sense considers a sequence of operators $-\ddiv(A_\eps \nabla \cdot)$ indexed by $\eps > 0$ and aims to characterize the limit as $\eps$ tends to zero. In many realistic applications, such a sequence of models can hardly be identified or may not be available in the first place. 
Assuming that the necessary requirements on the coefficient $A$ are met, a homogenized coefficient $A_\mathrm{hom}$ exists and does not involve oscillations on a fine scale. The coefficient $A_\mathrm{hom}$ can then be used in combination with a classical finite element ansatz, since $A_\mathrm{hom}$ does no longer include troublesome fine-scale quantities.
In practice, the homogenized coefficients cannot be computed easily and need to be approximated. This is, for instance, done with the Heterogeneous Multiscale Method (HMM)~\cite{EE03,EE05,AbdEEV12}, which in the end replaces $A$ by a computable approximation $A_h \in L^\infty(D, \mathbb{R}^{d\times d}_{\text{sym}})$ of $A_\mathrm{hom}$ with $(A_h)\vert_T \in \R^{d\times d}_{\text{sym}}$ for all $T \in \calT_h$. With this piecewise constant approximation of $A$, we obtain a possible compression operator $\mathfrak{C}$. Given an enumeration $1,\dots,m$ of the inner nodes in $\calT_h$ and writing $\lambda_1, \dots, \lambda_m$ for the associated nodal basis of $\Vh$, the compressed operator $\mathfrak{C}(A)$ can be defined as
\begin{equation}\label{eq:CHMM}
(\mathfrak{C}(A))_{i,j} = (S_{A})_{i,j} := \sum_{T \in \calT_h}\int_T ((A_h)\vert_T \nabla \lambda_j) \cdot \nabla \lambda_i \dx.
\end{equation}
That is, one takes the classical finite element stiffness matrix corresponding to the homogenized coefficient $A_\mathrm{hom}$ as an effective system matrix. In this case, the decomposition~\eqref{eq:decomp1} corresponds to a partition into element-wise stiffness matrices (with constant coefficient, respectively) that are merged with a simple finite element assembly routine.

 We emphasize that approaches based on analytical homogenization -- such as~\eqref{eq:CHMM} -- are able to provide reasonable approximations on the target scale $h$ but are subject to structural assumptions, in particular scale separation and local periodicity. The goal to overcome these restrictions has led to a new class of numerical methods that are specifically tailored to {treating} general coefficients with minimal assumptions. These methods are known as \emph{numerical homogenization approaches} and typically only require a boundedness condition as in \eqref{eq:frakAelliptic}.

\subsubsection{Compression by numerical homogenization}\label{sss:numhom}
The general idea of numerical homogenization methods is to replace the trial space $\Vh$ by a suitable \emph{multiscale space} $\tVh$, see for instance the references~\cite{MalP14,Owh17,OwhZB13,BabL11,AltHP21,Hou1997,EfeH09}. One possible construction uses a one-to-one correspondence of $\tVh$ to the space $\Vh$, which implies that the two spaces possess the same number of degrees of freedom. Typically, the multiscale space is chosen in a problem-adapted way. We indicate this dependence by defining the new space \mbox{$\tVh :=\calP_A \Vh$}, where $\calP_A\colon\Vh \to \V$ particularly depends on $A$.
Therefore, another possible choice of the operator $\mathfrak{C}$ leads to the effective matrix $\mathfrak{C}(A)$ given by
\begin{equation}\label{eq:CLOD}
(\mathfrak{C}(A))_{i,j} = (S_{A})_{i,j} := a_A(\calP_A \lambda_j, \lambda_i).
\end{equation}
A prominent example for such an approach -- and, thus, the operator $\mathfrak{C}$ -- is the Petrov--Galerkin version of the Localized Orthogonal Decomposition (LOD) method which explicitly constructs a suitable operator~$\calP_A$. The LOD was introduced in~\cite{MalP14} and theoretically and practically works for very general coefficients. It has also been successfully applied to other problem classes, for instance wave propagation problems in the context of Helmholtz and Maxwell equations~\cite{P17,GalP15,GalHV18,RenHB19,MaiV20} or the wave equation~\cite{AbdH17,PetS17,MaiP19,GeeM21}, eigenvalue problems~\cite{MalP15,MalP17}, and in connection with time-dependent nonlinear Schr\"odinger equations~\cite{HenW20}. However, it requires a slight deviation from locality. That is, while the classical finite element method and the HMM result in a system matrix that only includes neighbor-to-neighbor communication between the degrees of freedom, the multiscale approach~\eqref{eq:CLOD} moderately increases this communication to effectively incorporate the fine-scale information in $A$ for a broader range of coefficients, which is a common property of modern homogenization techniques. As indicated in~\cite{CaiMP20}, this slightly increased communication indeed seems to be necessary to handle very general coefficients. 

Since we consider a class $\mathfrak{A}$ of arbitrarily rough coefficients, the compression operator~\eqref{eq:CLOD} corresponding to the operator $\calP_A$ constructed in the LOD method is a suitable choice for our discussion as well as for the numerical experiments of Section~\ref{sec:numexp}. In the following subsection, we therefore {have a closer look} into its construction and summarize some main results. \corr{Note that we restrict ourselves to an elliptic model problem with homogeneous Dirichlet boundary conditions, but the compression approach can generally be extended to more involved settings such as the ones mentioned above.}

\subsection{Localized Orthogonal Decomposition}\label{subsec:LOD}
The method is based on a projective quasi-interpolation operator $\calI_h\colon \V\to \Vh$ with the following approximation and stability properties: for an element $T \in \calT_h$, we require that 
\begin{align*}
\|h^{-1}(v-\calI_h v)\|_{L^2(T)} + \|\nabla \calI_h v\|_{L^2(T)} \leq C \|\nabla v\|_{L^2(\Nb(T))}
\end{align*}
for all $v\in \V$, where the constant $C$ is independent of $h$, and $\Nb(S) := \Nb^1(S)$ is the neighborhood (of order $1$) of $S \subseteq D$ defined by
\begin{equation*}
\Nb^1(S) := \bigcup \bigl\{\overline{K} \in \calT_h\,\vert\, \overline{S} \,\cap\, \overline{K}\neq \emptyset\bigl\}.
\end{equation*} 
For a particular choice of $\calI_h$, we refer to \cite{ErnG15}.

For a given $\calI_h$ with the above properties, we can define the so-called fine-scale space $\W$, which contains all functions that are not well captured by the finite element functions in $\Vh$. It is defined as the kernel of $\calI_h$ with respect to $\V$, i.e.,
\begin{equation*}
\W := \ker {\calI_h}\vert_{\V},
\end{equation*}
and its local version, for any $S \subseteq D$, is given by
\begin{equation*}
\W(S) := \{w \in \W\,\vert\,\mathrm{supp}(w) \subseteq S\}.
\end{equation*}
In order to incorporate fine-scale information contained in the coefficient $A$, the idea is now to compute coefficient-dependent local corrections of functions $v_h\in\Vh$. To this end, we define the neighborhood of order $\ell\in\N$ iteratively by \corr{$\Nb^\ell(S) := \Nb(\Nb^{\ell-1}(S)),\,\ell \geq 2$}. For any function $v_h \in \Vh$, its element corrector $\calQ^\ell_{A,T} v_h{\in \W(\Nb^\ell(T))},\,T \in \calT_h$, is defined by
\begin{equation}\label{eq:corT}
a_A(\calQ^\ell_{A,T} v_h, w) = \int_T A \nabla v_h \cdot \nabla w \dx\quad  \text{for all } w \in \W(\Nb^\ell(T)).
\end{equation}
Note that in an implementation, the element corrections $\calQ^\ell_{A,T}$ have to be computed on a sufficiently fine mesh that resolves the oscillations of the coefficient $A$. Since the algebraic realization of the correctors and guidelines for an efficient implementation of the LOD method in general are not within the scope of the article, we refer to~\cite{EngHMP16} for details. We emphasize that, by construction, the supports of the correctors $\calQ^\ell_{A,T} v_h$ are limited to $\Nb^\ell(T)$.
The global correction $\calQ^\ell_A\colon \Vh \to \W$ then consists of a summation of these local contributions and is given by
\begin{equation*}
\calQ^\ell_A := \sum_{T \in \calT_h}\calQ^\ell_{A,T}.
\end{equation*}
Note that the choice $\ell = \infty$ corresponds to a computation of the element corrections on the entire domain $D$ and leads to the orthogonality property
\begin{equation}\label{eq:ortho}
a_A((1-\calQ^\infty_A)v_h,w) = 0\quad \text{for all } w\in \W,
\end{equation} 
that defines an $a_A$-orthogonal splitting $\V = (1-\calQ^\infty_A)\Vh \oplus \W$. This particularly explains the name \emph{orthogonal decomposition}. The use of localized element corrections is motivated by the decay of the corrections $\calQ^\infty_{A,T}$ away from the element $T$. This is, for instance, shown in~\cite{HenP13,Pet16} (based on \cite{MalP14}) and reads
\begin{equation*}
\|\nabla(\calQ^\infty_A-\calQ^\ell_A)v_h\|_{\LL} \leq Ce^{-c_\mathrm{dec}\ell}\, \|\nabla v_h\|_{\LL}
\end{equation*}
with a constant $c_\mathrm{dec}$ which is independent of $h$ and $\ell$.

Motivated by the decomposition~\eqref{eq:ortho} and the localized approximations in~\eqref{eq:corT}, we choose $\calP_A := 1 - \calQ^\ell_A$ in~\eqref{eq:CLOD}. 
The space $\tVh := \calP_A\Vh = (1 - \calQ^\ell_A)\Vh$, which has the same number of degrees of freedom as $\Vh$, can then be used as ansatz space for the discretization of \eqref{eq:modelelliptic}. Note that the original LOD method introduced in~\cite{MalP14} considers a discretization where $\tVh$ is also used as test space. We, however, consider the Petrov--Galerkin variant of the method as analyzed in~\cite{ElfGH15}, that uses the classical finite element space $\Vh$ as test space instead, i.e., we seek $u_h \in \Vh$ such that 
\begin{equation}\label{eq:LODPG}
a_A((1-\calQ^\ell_A)u_h,v_h) = {\langle f_h, v_h\rangle} \quad\text{for all }v_h \in \Vh,
\end{equation}
where $f_h = \mathfrak{M}f \in V_h$ is again a suitable approximation of $f \in H^{-1}(D)$. 
This defines a compressed operator $\mathfrak{S}_A$ as in~\eqref{eq:abstractmodeldisc} that maps $u_h\in \Vh$ to $f_h \in V_h$. 
As it turns out, the Petrov--Galerkin formulation has some computational advantages over the classical method, in particular in terms of memory requirement. For details we again refer to~\cite{EngHMP16}. The theory in~\cite{ElfGH15} shows that the {approximation $u_h$} defined in~\eqref{eq:LODPG} is {first-order accurate in $L^2(D)$} provided that $\ell \gtrsim |\log h|$ and, additionally, $f \in L^2(D)$. More precisely, it holds 
\begin{equation*}
\begin{aligned}
\|(-\ddiv(A \nabla))^{-1} - \mathfrak{S}_A^{-1} \ \|_{L^2(D)\rightarrow L^2(D)} &= \sup_{f\in L^2(D)} \frac{\|(-\ddiv(A \nabla))^{-1}(f) - \mathfrak{S}_A^{-1}(\mathfrak{M}f) \ \|_{L^2(D)}}{\| f\|_{L^2(D)}} \\
&\lesssim h\ + e^{-c_\mathrm{dec}\ell},
\end{aligned}
\end{equation*}
where $\mathfrak{M}f$ denotes the $L^2$-projection of $f$ onto $\Vh$. 
Note that the methodology can actually be applied to more general settings beyond the elliptic case, see for instance~\cite{MP20} for an overview. 

\subsection{System matrix surrogate}\label{subsec:surrogate}
We now return to the discussion of the compression operator $\mathfrak{C}$ introduced in \eqref{eq:CLOD}, that maps coefficients $A\in\mathfrak{A}$ to the effective system matrices
\begin{equation*}
(S_A)_{i,j} := a_A((1-\calQ^\ell_A)\lambda_j,\lambda_i)
\end{equation*}
obtained from the Petrov--Galerkin LOD method.
Once $S_A$ has been computed for a given $A$, an approximation $u_h = \sum_{j=1}^m U_j \lambda_j$ can be computed by solving the following linear system for the coefficients $U=(U_1,\dots,U_m)^T$:
\begin{equation*}
S_A U = F,
\end{equation*}
where $F := (\langle \mathfrak{M}f,\lambda_1\rangle,\ldots,\langle \mathfrak{M}f,\lambda_m\rangle)^T$. {Since $u_h$ is equivalently characterized by the solution of~\eqref{eq:LODPG}, it} captures the effective behavior of the solution to the continuous problem~\eqref{eq:modelelliptic} on the target scale $h$ as discussed in the previous subsection.

The remainder of this section is dedicated to showing how the abstract decomposition~\eqref{eq:decomp1} translates to the present LOD setting and how it can be implemented in practice. Writing $\mathcal{N}(S)$ for the set of inner mesh nodes on some subdomain $S\subseteq D$ and denoting $N_S = \vert\, \mathcal{N}(S)\vert$, the effective system matrix $S_A$ can be decomposed as 
\begin{equation}\label{eq:decompLOD1}
S_{A} = \sum_{T\in \calT_h} \Phi_{T} (S_{A,T}),
\end{equation}
where the matrices $S_{A,T}$ are local system matrices of the form 
\begin{equation}\label{eq:decompLOD2}
(S_{A,T})_{i,j} = \int_{\Nb^\ell(T)}A\ \nabla (1-\calQ^\ell_{A,T}) \lambda_j \cdot \nabla \lambda_i \dx, \quad j\in \mathcal{N}(T),\ i\in \mathcal{N}(\Nb^\ell(T)),
\end{equation}
i.e., they correspond to the interaction of the localized ansatz functions $(1-\calQ^\ell_{A,T}) \lambda_j$ associated with the nodes of the element $T$ with the classical first order nodal basis functions whose supports overlap with the {\emph{element neighborhood}} $\Nb^\ell(T)$. This means that $S_{A,T}$ is a $N_{\Nb^\ell(T)} \times {N_T}$ 
matrix. In practice, the coefficient $A$ in~\eqref{eq:decompLOD2} is often replaced with an element-wise constant approximation $A_\eps$ on a finer mesh $\calT_\varepsilon$ that resolves all the oscillations of $A$ and that we assume to be a uniform refinement of $\calT_h$. 

As already explained in the abstract framework, the mappings $\Phi_{T}$ in \eqref{eq:decompLOD1} are local-to-global mappings that assemble the contributions $S_{A,T}$ {on an element neighborhood} to a global matrix. In particular, given an enumeration $1,\dots ,N_{\Nb^\ell(T)}$ of the nodes in $\mathcal{N}(\Nb^\ell(T))$ one considers a mapping $g_T(\cdot)$ that assigns to a given node index $i$ in the {element neighborhood} $\Nb^\ell(T)$ its global node index $g_{T}(i)\in\{1,\dots,m\}$. This mapping can be represented by an $m\times N_{\Nb^\ell(T)}$ sparse matrix with entries
\begin{equation*}
\pi_T[i,j] = \begin{cases} 1,\ &\text{if } i = g_{T}(j), \\0,\ &\text{otherwise.} \end{cases}
\end{equation*}
Analogously, given an enumeration $1,\dots,{N_T}$ 
of the nodes in $\mathcal{N}(T)$, there exists a mapping $\tilde{g}_T(\cdot)$ -- represented by an $m\times {N_T}$-matrix -- that assigns to a given node in $\mathcal{N}(T)$ with index $i$ its global representative with index $\tilde{g}_{T}(i)$. The corresponding matrix is given by
\begin{equation*}
\varphi_T[i,j] = \begin{cases} 1,\ &\text{if } i = \tilde{g}_{T}(j), \\0,\ &\text{otherwise.} \end{cases}
\end{equation*}
Using these matrices, the decomposition \eqref{eq:decompLOD1} reads
\corr{
\begin{equation}\label{eq:decompLOD3}
S_{A} = \sum_{T\in \calT_h} \pi_T\ S_{A,T}\ \varphi_T^{\, \prime},
\end{equation}
where $\varphi_T^{\, \prime}$ denotes the transpose of the matrix $\varphi_T$.
}

From the definition of the local contributions $S_{A,T}$ introduced in~\eqref{eq:decompLOD2}, it directly follows that $S_{A,T}$ does only depend on the restriction of $A$, respectively $A_\eps$, to the element neighborhood $\Nb^\ell(T)$. 
Let now $\calT_\varepsilon(\Nb^\ell(T))$ be the restriction of the mesh $\calT_\eps$ to $\Nb^\ell(T)$, consisting of $r = |\calT_\varepsilon(\Nb^\ell(T))|$ elements. Enumerating the elements then leads to the following operators that correspond to the abstract reduction operators in~\eqref{eq:reduction},
\begin{equation}\label{eq:reductionLOD}
R_{T}\colon \mathfrak{A} \rightarrow \R^{r},
\end{equation}
that map a global coefficient $A$ to a vector that contains the values of $A_\eps$ in the respective cells of $\calT_\varepsilon(\Nb^\ell(T))$. 

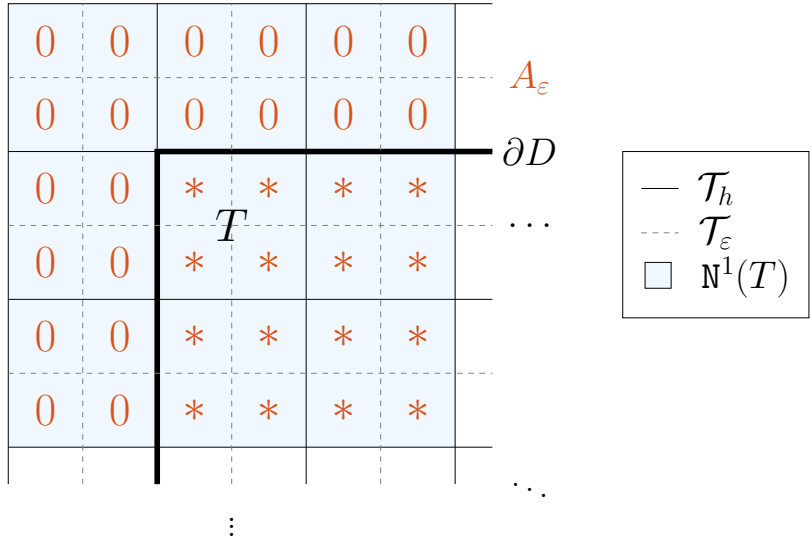
\begin{figure}
\scalebox{0.7}{\begin{tikzpicture}[scale=0.7]

\fill [aliceblue] (0,0) rectangle (12,12);

\draw (16.5,3.5) -- (16.5,8) -- (21.5,8) -- (21.5,3.5) -- cycle;
\draw (17, 7) -- (18,7);
\node[] at (19,7) {\huge{$\mathcal{T}_h$}}; 
\draw[gray, dashed] (17, 5.8) -- (18,5.8);
\node[] at (19,5.8) {\huge{$\mathcal{T}_\varepsilon$}}; 
\fill[aliceblue] (17.1,4.3) rectangle (17.8,5);
\draw[] (17.1,4.3) rectangle (17.8,5);
\node[] at (19.8,4.6) {\huge{$\Nb^1(T)$}};

\draw (0,-1) -- (0,12) -- (13,12);

\draw[line width = 3] (4,-1) -- (4,8) -- (13,8);
\node[] at (14, 8) {\huge{$\partial D$}};

\node[] at (14, 6) {\huge{$\hdots$}};
\node[] at (6, -2) {\huge{$\vdots$}};
\node[] at (14, -1) {\huge{$\ddots$}};

\draw (4,-1) -- (4,12);
\draw (8,-1) -- (8,12);
\draw (12,-1) -- (12,12);
\draw (0,0) -- (13,0);
\draw (0,4) -- (13,4);
\draw (0,8) -- (13,8);

\draw[gray, dashed] (2,-1) -- (2,12);
\draw[gray, dashed] (6,-1) -- (6,12);
\draw[gray, dashed] (10,-1) -- (10,12);
\draw[gray, dashed] (0,2) -- (13,2);
\draw[gray, dashed] (0,6) -- (13,6);
\draw[gray, dashed] (0,10) -- (13,10);

\node[] at (6,6) {\Huge{\textbf{$T$}}};

\node[myOrange] at (14, 10) {\huge{$A_\eps$}};
\node[myOrange] at (1,1) {\Huge{\textbf{$0$}}};
\node[myOrange] at (1,3) {\Huge{\textbf{$0$}}};
\node[myOrange] at (1,5) {\Huge{\textbf{$0$}}};
\node[myOrange] at (1,7) {\Huge{\textbf{$0$}}};
\node[myOrange] at (1,9) {\Huge{\textbf{$0$}}};
\node[myOrange] at (1,11) {\Huge{\textbf{$0$}}};
\node[myOrange] at (3,1) {\Huge{\textbf{$0$}}};
\node[myOrange] at (3,3) {\Huge{\textbf{$0$}}};
\node[myOrange] at (3,5) {\Huge{\textbf{$0$}}};
\node[myOrange] at (3,7) {\Huge{\textbf{$0$}}};
\node[myOrange] at (3,9) {\Huge{\textbf{$0$}}};
\node[myOrange] at (3,11) {\Huge{\textbf{$0$}}};
\node[myOrange] at (5,11) {\Huge{\textbf{$0$}}};
\node[myOrange] at (7,11) {\Huge{\textbf{$0$}}};
\node[myOrange] at (9,11) {\Huge{\textbf{$0$}}};
\node[myOrange] at (11,11) {\Huge{\textbf{$0$}}};
\node[myOrange] at (5,9) {\Huge{\textbf{$0$}}};
\node[myOrange] at (7,9) {\Huge{\textbf{$0$}}};
\node[myOrange] at (9,9) {\Huge{\textbf{$0$}}};
\node[myOrange] at (11,9) {\Huge{\textbf{$0$}}};

\node[myOrange] at (5,1) {\Huge{\textbf{$\ast$}}};
\node[myOrange] at (5,3) {\Huge{\textbf{$\ast$}}};
\node[myOrange] at (5,5) {\Huge{\textbf{$\ast$}}};
\node[myOrange] at (5,7) {\Huge{\textbf{$\ast$}}};
\node[myOrange] at (7,1) {\Huge{\textbf{$\ast$}}};
\node[myOrange] at (7,3) {\Huge{\textbf{$\ast$}}};
\node[myOrange] at (7,5) {\Huge{\textbf{$\ast$}}};
\node[myOrange] at (7,7) {\Huge{\textbf{$\ast$}}};
\node[myOrange] at (9,1) {\Huge{\textbf{$\ast$}}};
\node[myOrange] at (9,3) {\Huge{\textbf{$\ast$}}};
\node[myOrange] at (9,5) {\Huge{\textbf{$\ast$}}};
\node[myOrange] at (9,7) {\Huge{\textbf{$\ast$}}};
\node[myOrange] at (11,1) {\Huge{\textbf{$\ast$}}};
\node[myOrange] at (11,3) {\Huge{\textbf{$\ast$}}};
\node[myOrange] at (11,5) {\Huge{\textbf{$\ast$}}};
\node[myOrange] at (11,7) {\Huge{\textbf{$\ast$}}};

\end{tikzpicture} }
\caption{Illustration of the extended {element neighborhood} $\Nb^1(T)$ around a corner element $T\in\calT_h$. An asterisk indicates that $A_\eps {\vert_K} \in [\alpha,\beta]$, a zero that $A_\eps{\vert_K}  = 0$ in the respective cell {$K$} of the refined mesh $\calT_\varepsilon$.}\label{fig:boundarypatch}
\end{figure}

As already mentioned in the abstract section above, we aim for a uniform output size of the operators $R_T$, since the outputs of the operators $R_T$ will later on be fed into a neural network with a fixed number of input neurons. In order to achieve that, we artificially extend the {domain $D$ and the} mesh $\calT_h$ by $\ell$ layers of outer elements around the boundary elements of $\calT_h$, thus ensuring that the {element neighborhood} $\Nb^\ell(T)$ always consists of the same number of elements regardless of the respective location of the central element $T\in\calT_h$ relative to the boundary. Further, we extend the piecewise constant coefficient $A_\eps$ by zero on those outer elements. Figure~\ref{fig:boundarypatch} illustrates this procedure for the case $d=2$ and $\ell=1$ for an element $T$ that lies in a corner of the computational domain. In this figure, $\calT_h$ is a uniform quadrilateral mesh on the domain $D$ and $\calT_\varepsilon$ is obtained from $\calT_h$ by a single uniform refinement step.
The asterisks indicate the coefficient $A_\eps$ taking a regular value in the interval $[\alpha,\beta]$, whereas in the cells outside of $D$, we set $A_\eps$ to zero.

Note that this enlargement of the mesh $\calT_h$ to obtain equally sized element neighborhoods $\Nb^\ell(T)$ also introduces artificial mesh nodes that lie outside of $D$ {and that are all formally considered as \emph{inner} nodes for the definition of $N_S = \vert\, \mathcal{N}(S)\vert$ with a subset $S$ in the extended domain}. This implies that the local system matrices $S_{A,T}$ of dimension $N_{\Nb^\ell(T)} \times {N_T}$ introduced in \eqref{eq:decompLOD2} are all of equal size as well and the rows of $S_{A,T}$ corresponding to test functions associated with nodes {that are attached to outer elements} contain only zeros. During the assembly process of the local contributions to a global matrix, these zero rows are disregarded (which is also consistent with our definition of the matrices $\pi_T,\ \varphi_T$). 

Finally, in order to unify the computation of local contributions, we use an abstract mapping $\mathfrak{C}_{\text{red}}$ with fixed input dimension $r$ and fixed output dimension $N_{\Nb^\ell(T)} \times 2^d$ as proposed for the abstract framework in Section~\ref{subsec:matrixdecomp}. The mapping takes the restriction of $A_\eps$ to an {element neighborhood} $\Nb^\ell(T)$ as input data and outputs the corresponding approximation of a local effective matrix $S_{A,T}$ that will be determined by an underlying neural network $\Psi(\cdot,\theta)$.

%
%

\section{Numerical Experiments}\label{sec:numexp}
In this section, we demonstrate the feasibility of our proposed approach by performing numerical experiments in the setting of Section \ref{sec:numhom}. For all experiments, we consider the two-dimensional computational domain $D = (0,1)^2$, which we discretize with a uniform quadrilateral mesh $\calT_h$ of characteristic mesh size $h = 2^{-5}$. The coefficients are allowed to vary on the finer unresolved scale $\varepsilon = 2^{-8}$.

\subsection{Coefficient family} \label{subsec: coefclass}
In order to test our method's ability to deal with coefficients that show oscillating behavior across multiple scales, we introduce a hierarchy of meshes $\calT^{k},\ k = 0, 1,\dots,8,$ where the initial mesh $\calT^0$ consists only of a single element, and the subsequent meshes are obtained by uniform refinement, i.e., $\calT^k$ is obtained from $\calT^{k-1}$  by subdiving each element of $\calT^{k-1}$ into four equally sized elements. This implies that the characteristic mesh size of $\calT^k$ is given by $2^{-k}$. In the following, we refer to the parameter $k$ as the \emph{mesh level}. In particular, the computational mesh $\calT_h = \calT^5$ corresponds to the mesh level $5$, whereas the coefficients may vary on the mesh level $8$ and are therefore only resolved by the finest mesh $\calT^8$. We thus have a scenario where an information gap of $3$ mesh levels has to be bridged. Based on the mesh hierarchy, we now define the coefficient family $\mathfrak{A}$ of interest. Let $\mathfrak{A}_k$ denote the set of element-wise constant coefficients on $\calT^k$ whose values in the respective cells are iid uniformly distributed on the interval $[\alpha,\beta]:=[1,5]$, i.e., 
\begin{equation*}
\mathfrak{A}_k := \left\{ A\ \in Q^0(\calT^k)\ \big\vert\ A\vert_{T} \overset{\text{iid}}{\sim} U([1,5])\ \text{for all}\ T\in\calT^k \right\}.
\end{equation*}
Furthermore, let $\mathfrak{A}_{\text{ms}}$ denote the set of coefficients of the form
\begin{equation*}
A = \frac{1}{9} \sum_{k=0}^8 A_k,\ A_k\in \mathfrak{A}_k.
\end{equation*}
These multiscale coefficients are especially interesting, since they vary on all considered scales simultaneously and are therefore very hard to deal with using classical homogenization techniques due to their unstructured nature. The total set of interest $\mathfrak{A}$ is then defined as
\begin{equation*}
\mathfrak{A} := \bigcup_{k = 0}^8 \mathfrak{A}_k \cup \mathfrak{A}_{\text{ms}}.
\end{equation*}
In the following, we will frequently index coefficients sampled from $\mathfrak{A}$ by their corresponding level, i.e., write $A_k,\ k\in\{0,\dots,8,\text{ms}\}$ instead of a plain $A$.

\subsection{Data generation and preprocessing}\label{subsec:preproc}
In order to train the network, we sample $500$ coefficients $A_k^{(i)},\ i = 1,\dots,500$ from each $\mathfrak{A}_k,\ k\in\{0,\dots,8,\text{ms}\}$, where the individual samples $A_k^{(i)}$ on the coarser mesh levels $k = 0,\dots,7$ are prolongated to the finest mesh $\calT^8$ in order to achieve a uniform dimension across all scales. The set of all sample coefficients is subsequently divided into a training, validation, and test set according to a $80-10-10$ split. In order to achieve an identical distribution in all three sets, the splitting is performed separately on every level (including ms), i.e., for every $k\in \{0,\dots,8,\text{ms}\}$, the first $400$ coefficients $A_k^{(i)}, i = 1,\dots,400$ get assigned to the training set $\calD_{\text{train}}$, the samples with indices $401, \dots, 450$ to the validation set $\calD_{\text{val}}$ and those with indices $451, \dots, 500$ to the test set $\calD_{\text{test}}$. Then, we individually split each sample $A_k^{(i)},$ using the reduction operators $R_T$ introduced in \eqref{eq:reductionLOD}, into sub-samples $A_{k,T}^{(i)}$ based on {element neighborhoods} $\Nb^\ell(T)$ for $\ell = 2$ that are centered around the elements $T\in\calT_h$, also taking into account the artificial extension of the {element neighborhoods} around the boundary of $D$. Since our target scale of interest is $h = 2^{-5}$ and $\calT_h$ is a uniform quadrilateral mesh, this yields $1024$ sub-samples $A_{k,T}^{(i)} \in \mathbb{R}^{1600}$ per sample $A_k^{(i)} \in \mathbb{R}^{65536}$. \corr{Note that the size of the sub-samples is obtained from the construction of the local neighborhoods $\Nb^\ell(T)$. Here, each neighborhood consists of $(2\ell + 1)^2 = 25$ elements in $\calT_h = \calT^5$, which corresponds to $64\cdot 25 = 1600$ elements in the mesh $\calT_\eps = \calT^8$.}

The corresponding ``labels'', i.e., the local effective system matrices $S_{A,k,T}^{(i)} \in \mathbb{R}^{36 \times 4}$, are then computed with the Petrov--Galerkin LOD according to~\eqref{eq:decompLOD2} and flattened column-wise to vectors in $\mathbb{R}^{144}$. In total, we obtain $10 \cdot 400 \cdot 1024 = 4096000$ pairs $(A_{k,T}^{(i)}, S_{A,k,T}^{(i)}) \in \calD_{\text{train}}$ to train our network with, and $512000$ pairs in $\calD_{\text{val}}$ and $\calD_{\text{test}}$ each.

\subsection{Network architecture and training}\label{subsec:network}
Given the above dataset, we now try to fit it with a suitable neural network $\Psi$. As network architecture, we consider a dense feedforward network with a total of eight layers including the input and output layer. As activation function, we choose the standard ReLU activation given by $\rho(x):= \max(0,x)$ in the first seven layers and the identity function in the last layer. By convention, the activation function acts component-wise on vectors. The network output is thus of the form
\begin{equation}\label{eq:architecture1}
\Psi(x) = W^{(8)}\rho(W^{(7)}(\dots \rho(W^{(2)}\rho(W^{(1)}x + b^{(1)})+b^{(2)})\dots) + b^{(7)}) + b^{(8)} ,
\end{equation}
where the weight matrices and bias vectors have the following dimensions:
\begin{equation*}\label{eq:architecture2}
\begin{aligned}
&W^{(1)} \in \mathbb{R}^{1600\times 1600}, \quad &&W^{(2)} \in \mathbb{R}^{800\times 1600}, \quad
&&W^{(3)} \in \mathbb{R}^{800\times 800}, \quad &&W^{(4)} \in \mathbb{R}^{400\times 800}, \\ 
&W^{(5)} \in \mathbb{R}^{400\times 400}, \quad &&W^{(6)} \in \mathbb{R}^{144\times 400}, \quad
&&W^{(7)} \in \mathbb{R}^{144\times 144}, \quad &&W^{(8)} \in \mathbb{R}^{144\times 144}, \\
&b^{(1)} \in \mathbb{R}^{1600}, \quad &&b^{(2)}, \in \mathbb{R}^{800}, \quad
&&b^{(3)} \in \mathbb{R}^{800}, \quad &&b^{(4)} \in \mathbb{R}^{400}, \\ 
&b^{(5)} \in \mathbb{R}^{400}, \quad &&b^{(6)} \in \mathbb{R}^{144}, \quad
&&b^{(7)} \in \mathbb{R}^{144}, \quad &&b^{(8)} \in \mathbb{R}^{144},
\end{aligned}
\end{equation*}
\begin{figure}
     \centering
     \includegraphics[width=0.5\textwidth]{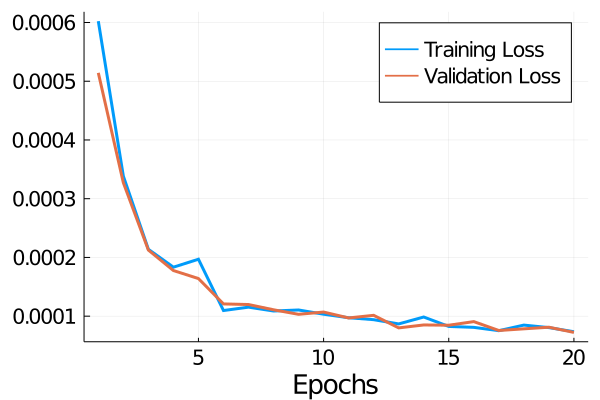}
\caption{Development of loss functional $\mathcal{J}$ during $20$ epochs of the network training.}
\label{fig:loss}
\end{figure}
yielding a total of $5063504$ trainable parameters. The idea behind this architecture is that in the first six layers, information about the coefficient in the input {element neighborhood} is gathered and combined by allowing communication between all inputs in layers with odd indices, whereas in the layers with even indices, this information is repeatedly compressed. \corr{That is, every other layer is built in such a way that the input and output dimension are equal. If the neurons in that layer are understood as some sort of degrees of freedom in a mesh, this refers to having communication among all of these degrees of freedoms, while the layers in between reduce the number of degrees of freedom, which can be interpreted as transferring information to a coarser mesh.} In the last two layers, this compressed information is taken and assembled to the local effective system matrix. Note that this logarithmic dependence of the number of layers on the number of scales that need to be bridged by the network (two layers per mesh level to be bridged plus two layers to assemble the local effective matrix) yielded the most reliable results in our experiments. Shallower networks had difficulties fitting the complex training set consisting of coefficients varying on different scales, whereas deeper networks were more prone to overfitting the training set. More involved architectures, for example ones that include skip connections between layers like in the classic ResNet~\cite{HeZRS16} are also conceivable, however, this seems not to be necessary to obtain good results. \corr{The key message here is that the coefficient-to-surrogate map can be satisfyingly approximated by a simple feedforward architecture, whose size does depend only on the scales $\varepsilon$ and $h$, but not on any finer discretization scales.} The implementation of the network as well as the training is performed using the library \emph{Flux}~\cite{Inn18} for the open-source scientific computing language \emph{Julia}~\cite{BezEKS17}. 

\begin{figure}
     \centering
     \begin{subfigure}[a]{0.49\textwidth}
         \centering
         \includegraphics[width=\textwidth]{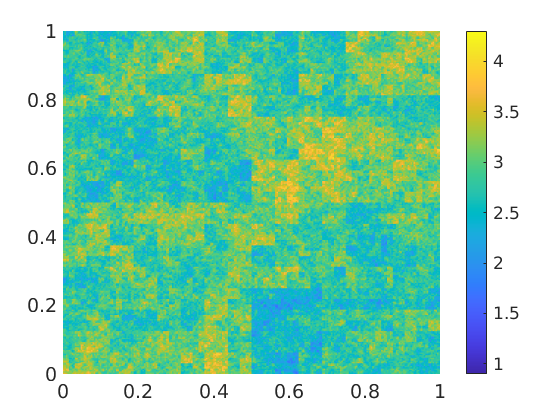}
     \end{subfigure} 
          \begin{subfigure}[a]{0.49\textwidth}
         \centering
         \includegraphics[width=\textwidth]{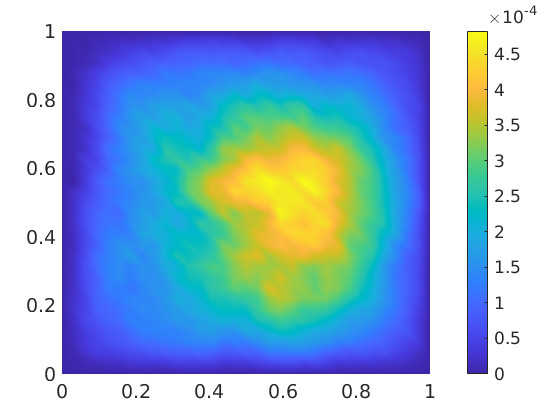}
     \end{subfigure} 
     \begin{subfigure}[a]{0.49\textwidth}
         \centering
         \includegraphics[width=\textwidth]{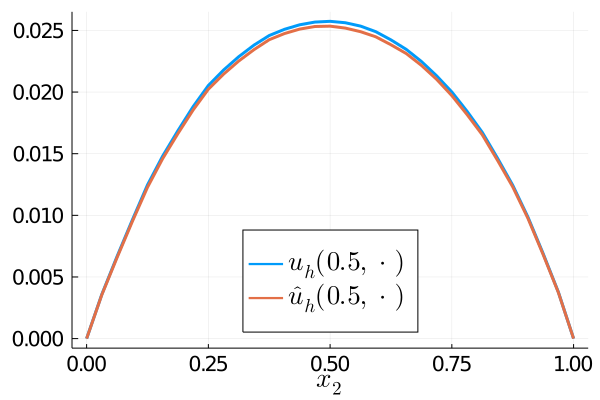}
     \end{subfigure}
     \begin{subfigure}[a]{0.49\textwidth}
         \centering
         \includegraphics[width=\textwidth]{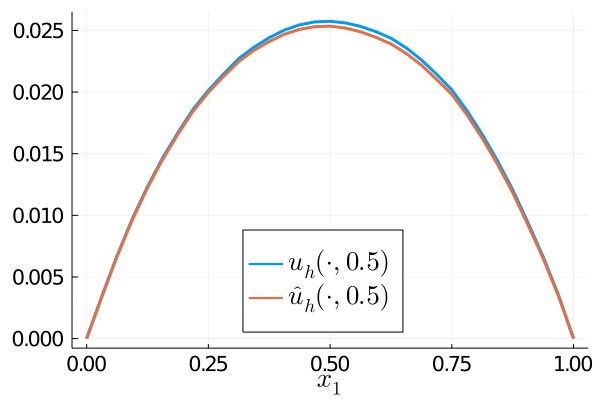}
     \end{subfigure}  
\caption{Results of Experiment 1: unstructured multiscale coefficient sampled from $\mathfrak{A}_{\text{ms}}$ (top left), $\vert u_h(x) - \widehat{u}_h(x) \vert$ (top right) and comparison of $u_h$ vs. $\widehat{u}_h$ along the cross sections $x_1 = 0.5$ (bottom left) and $x_2 = 0.5$ (bottom right).}
\label{fig:expmulti}
\end{figure}
\begin{figure}
     \centering
     \begin{subfigure}[a]{0.49\textwidth}
         \centering
         \includegraphics[width=\textwidth]{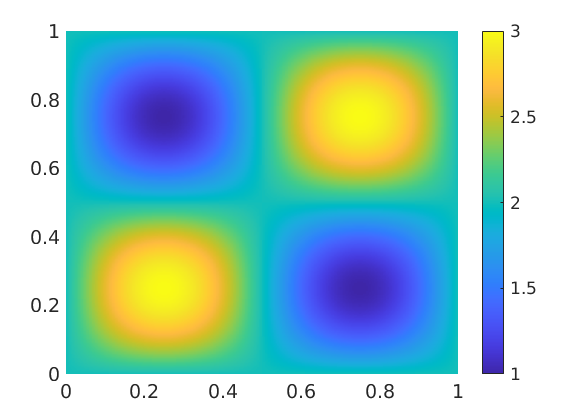}
     \end{subfigure} 
          \begin{subfigure}[a]{0.49\textwidth}
         \centering
         \includegraphics[width=\textwidth]{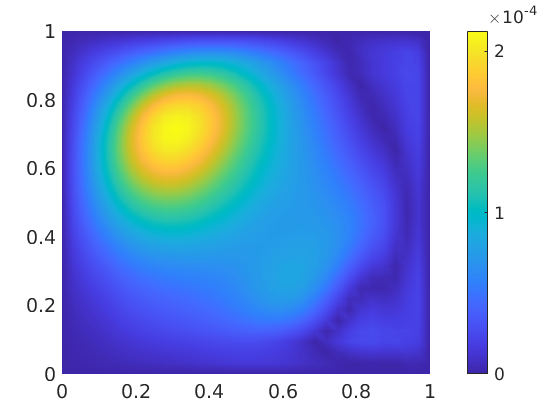}
     \end{subfigure} 
     \begin{subfigure}[a]{0.49\textwidth}
         \centering
         \includegraphics[width=\textwidth]{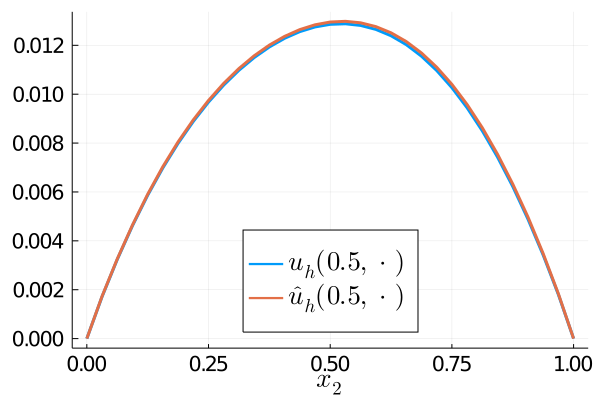}
     \end{subfigure}
     \begin{subfigure}[a]{0.49\textwidth}
         \centering
         \includegraphics[width=\textwidth]{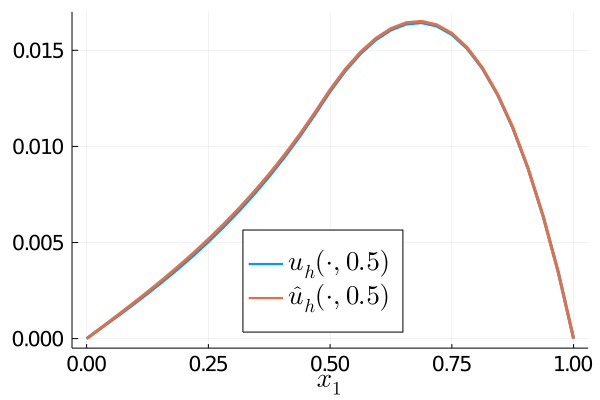}
     \end{subfigure}  
\caption{Results of Experiment 2: smooth coefficient (top left), $\vert u_h(x) - \widehat{u}_h(x) \vert$ (top right) and comparison of $u_h$ vs. $\widehat{u}_h$ along the cross sections $x_1 = 0.5$ (bottom left) and $x_2 = 0.5$ (bottom right).}
\label{fig:expsine}
\end{figure}
\begin{figure}
     \centering
     \begin{subfigure}[a]{0.49\textwidth}
         \centering
         \includegraphics[width=\textwidth]{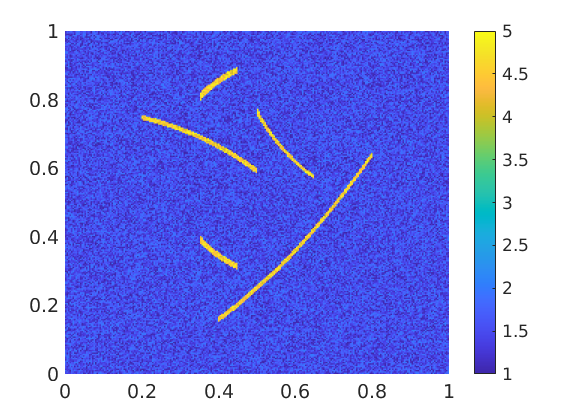}
     \end{subfigure} 
          \begin{subfigure}[a]{0.49\textwidth}
         \centering
         \includegraphics[width=\textwidth]{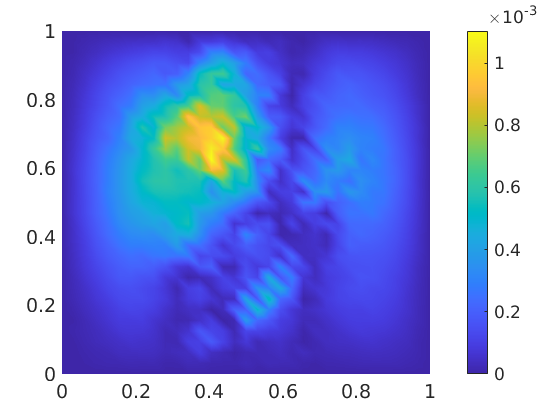}
     \end{subfigure} 
     \begin{subfigure}[a]{0.49\textwidth}
         \centering
         \includegraphics[width=\textwidth]{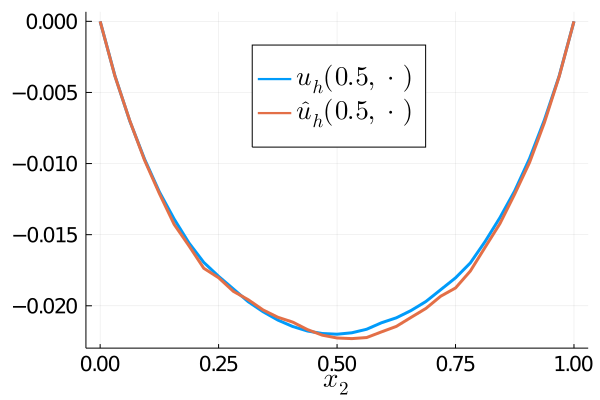}
     \end{subfigure}
     \begin{subfigure}[a]{0.49\textwidth}
         \centering
         \includegraphics[width=\textwidth]{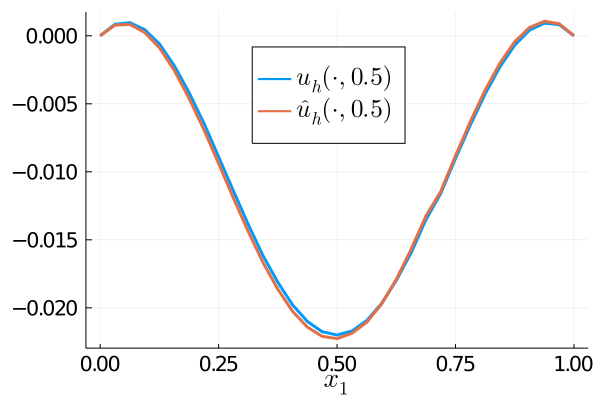}
     \end{subfigure}  
\caption{Results of Experiment 3: coefficient with cracks (top left),  $\vert u_h(x) - \widehat{u}_h(x) \vert$ (top right) and comparison of $u_h$ vs. $\widehat{u}_h$ along the cross sections $x_1 = 0.5$ (bottom left) and $x_2 = 0.5$ (bottom right).}
\label{fig:expcrack}
\end{figure}

After initializing all parameters in the network according to a Glorot uniform~ distribution~\cite{GloB10}, the network \eqref{eq:architecture1} is trained on minibatches of $1000$ samples for a total of $20$ epochs on $\calD_{\text{train}}$, using the ADAM optimizer~\cite{KinB14} with a step size of $10^{-4}$ for the first $5$ epochs before reducing it to $10^{-5}$ for the subsequent $15$ epochs. \corr{It could be observed that further training led to a stagnation of the validation error, whereas the error on the training set continued to decrease (very slowly but gradually), indicating overfitting of the network.} The development of the loss functional $\mathcal{J}$ defined in \eqref{eq:loss} over the epochs is shown in Figure \ref{fig:loss}. Note that training and validation loss stay very close to each other during the whole training process since $\calD_{\text{train}}$ and $\calD_{\text{val}}$ have the same sample distribution due to our chosen splitting procedure.

The development of the loss during the training and an average loss of $7.78\cdot 10^{-5}$ on the test set $\calD_{\text{test}}$ indicates that the network has at least learned to approximate the local effective system matrices. In applications, however, we are mostly concerned about how well this translates to the global level, when computing solutions to problem \eqref{eq:modelelliptic} using a global system matrix assembled from network outputs. In order to investigate this question, the next three subsections are dedicated to evaluating the performance of the trained network at exactly this task for several coefficients unseen during training. For a given right-hand side $f$ and coefficient $A$, we denote with $u_h$ the solution of~\eqref{eq:LODPG}, obtained with the Petrov-Galerkin LOD matrix $S_A$ defined in \eqref{eq:decompLOD3}, and with $\widehat{u}_h$ the solution obtained by using the neural network approximation of this matrix, i.e.,
\corr{
\begin{equation}\label{eq:ShatNN}
\widehat{S}_{A} := \sum_{T\in \calT_h} \pi_T\ \Psi(R_T(A))\ \varphi_T^{\, \prime}.
\end{equation}
}

\corr{The spectral norm difference $\|S_A - \widehat{S}_A\|_2$,} the $L^2$-error $\|u_h - \widehat{u}_h \|_{\LL}$ as well as the visual discrepancy between the two solutions are then considered as a measure of the network's global performance. \corr{We emphasize once more that computing approximate surrogates via \eqref{eq:ShatNN} is significantly faster compared to~\eqref{eq:decompLOD1} and~\eqref{eq:decompLOD2}. This is due to the fact that no corrector problems of the form \eqref{eq:corT} have to be solved to obtain the surrogate model. As pointed out, the solution of these local problems requires inversion on a very fine discretization scale that is significantly smaller than the scale $\varepsilon$ on which the coefficient varies. In order to compute the system matrix~$S_A,$ one has to solve $N_{\mathcal{T}_h}$ fine-scale linear systems, where $N_{\mathcal{T}_h}$ denotes the number of elements in $\mathcal{T}_h$. In contrast, the main computational effort of evaluating our trained network consists of $N_L$ matrix-matrix multiplications, where $N_L$ is the number of layers in the network (not taking into account bias vectors and the activation function).}

\subsection{Experiment 1}\label{subsec:exp1}
For our first experiment we consider an unstructured multiscale coefficient sampled from $\mathfrak{A}_{\text{ms}}$ that was not a part of the training set and a constant right-hand side $f\equiv 1$. The coefficient (top left), the error $\vert u_h(x) - \widehat{u}_h(x)\vert$ (top right) as well as representative cross-sections along $x_1 = 0.5$ (bottom left) and $x_2 = 0.5$ (bottom right) of the two solutions $u_h$ and $\widehat{u}_h$ are shown in Figure \ref{fig:expmulti}. \corr{The spectral norm difference $\|S_A - \widehat{S}_A\|_2 \approx 6.58\cdot 10^{-2}$} and the $L^2$-error \corr{$\|u_h - \widehat{u}_h \|_{\LL} \approx 2.13 \cdot 10^{-4}$} confirm the visual impression -- the network has successfully learned to produce a system matrix that is able to capture the behavior of the solution on the target scale well.

\subsection{Experiment 2}\label{subsec:exp2}
Next, we test the network's performance for smoother and more regular coefficients than the ones it has been trained with. As a demonstrating example, we consider the coefficient $A(x) = 2 + \sin(2\pi x_1) \sin(2\pi x_2)$. The network's input is obtained by evaluating $A$ on the midpoints of the mesh $\calT_\varepsilon$ on the fine unresolved scale $\eps$. In this example, we choose the function $f(x) = x_1 \chi_{\{x_1 \geq 0.5 \}}$ as right-hand side, where {$\chi_S$} denotes the characteristic function of {the set $S \subseteq D$}. The results are shown in Figure \ref{fig:expsine}. We obtain an even better $L^2$-error of \corr{$\|u_h - \widehat{u}_h \|_{\LL} \approx 7.56 \cdot 10^{-5}$} and \corr{a spectral norm difference of $\|S_A - \widehat{S}_A\|_2 \approx 3.91\cdot 10^{-2}$.} A comparison between the solutions at the cross-sections shows that there is almost no discernible visual difference between the LOD-solution and the approximation obtained using our trained network. 

\subsection{Experiment 3}\label{subsec:exp3}
As a third experiment, we choose another coefficient that possesses an unfamiliar structure not seen by the network during the training phase, this time a less regular one. The coefficient is shown in the top left of Figure \ref{fig:expcrack}. It is composed of a background part (blue region), which is obtained by sampling uniformly and independently on each element of $\calT_\varepsilon$ from the interval $[1,2]$, and several ``cracks'' (yellow regions), in which the coefficient varies uniformly in $[4,5]$. The right-hand side here is $f(x) = \cos(2\pi x_1)$. A computation of the $L^2$-error \corr{$\|u_h - \widehat{u}_h \|_{\LL} \approx 2.72 \cdot 10^{-4}$} shows that the overall error is still moderate, an a closer visual inspection of the solutions along the cross-sections however reveals more prominent deviations of the neural network approximation to the ground truth. \corr{The spectral norm difference $\|S_A - \widehat{S}_A\|_2 \approx 2.81\cdot 10^{-1}$ is also one order of magnitude larger than in the previous examples.} \corr{Nevertheless, it might be possible that performing a few corrective training steps including samples of this nature would be sufficient to fix this issue. A thorough investigation of this hypothesis, along with the extension to other coefficient classes is subject of future work.}

%
%

\section{Conclusion and Outlook}\label{sec:conclusion}
We proposed an approach to the compression of linear differential operators -- parameterized by PDE coefficients that may depend on microscopic quantities that are not resolved by a target discretization scale of interest -- to lower-dimensional surrogates that are based on a combination of existing model reduction methods with a data-driven deep learning framework. Our method is motivated by the fact that the computation of the surrogates (represented by effective system matrices) via classical methods is nontrivial and requires significant computational resources in multi-query settings. To overcome this problem, we showed how the compression process can be approximated by a neural network based on given training data that can be generated using existing compression approaches. {Importantly, we avoid a global approximation by a neural network and instead first decompose the compression map into local contributions, which can then be approximated by one single unified network.} As an example, we studied a class of second-order elliptic diffusion operators. We showed how to approximate the compression map based on the Petrov--Galerkin formulation of the Localized Orthogonal Decomposition method with a neural network. The proposed ansatz has been numerically validated for a large set of piecewise constant and highly oscillatory multiscale coefficients. Furthermore, it has been shown that the approach also generalizes well, in the sense that a well-trained network is able to produce reasonable results even for classes of coefficients that it has not been trained on.

For future research, many possible research questions building on the present work are conceivable. Straightforward extensions would be to consider stochastic settings with differential operators parameterized by random fields or settings with high contrast. Another question to investigate is to what degree the method can be made robust against changes in geometry, for example by training the network not only on coefficients that are sampled on a fixed domain, but rather on reference patches with varying geometries. Mimicking a hierarchical discretization approach, one may also try to directly approximate the inverse operator which can be represented by a sparse matrix~\cite{FeiP20}. On a more theoretical level, the approximation properties of neural networks for various existing compression operators could be investigated, along with the question of the number of training samples required to faithfully approximate those for a given family of coefficients.

%
%

\newcommand{\etalchar}[1]{$^{#1}$}


\begin{thebibliography}{WCC{\etalchar{+}}20}

\bibitem[AAB{\etalchar{+}}15]{AbaAB+15}
M.~Abadi, A.~Agarwal, P.~Barham, E.~Brevdo, Z.~Chen, C.~Citro, G.~S. Corrado,
  A.~Davis, J.~Dean, M.~Devin, S.~Ghemawat, I.~Goodfellow, A.~Harp, G.~Irving,
  M.~Isard, Y.~Jia, R.~Jozefowicz, L.~Kaiser, M.~Kudlur, J.~Levenberg,
  D.~Man\'{e}, R.~Monga, S.~Moore, D.~Murray, C.~Olah, M.~Schuster, J.~Shlens,
  B.~Steiner, I.~Sutskever, K.~Talwar, P.~Tucker, V.~Vanhoucke, V.~Vasudevan,
  F.~Vi\'{e}gas, O.~Vinyals, P.~Warden, M.~Wattenberg, M.~Wicke, Y.~Yu, and
  X.~Zheng.
\newblock {TensorFlow}: Large-scale machine learning on heterogeneous systems,
  2015.
\newblock Software available from tensorflow.org.

\bibitem[ABS{\etalchar{+}}20]{ArbBSRK20}
H.~Arbabi, J.~E. Bunder, G.~Samaey, A.~J. Roberts, and I.~G. Kevrekidis.
\newblock Linking machine learning with multiscale numerics: Data-driven
  discovery of homogenized equations.
\newblock {\em Jom}, 72(12):4444--4457, 2020.

\bibitem[AEEV12]{AbdEEV12}
A.~Abdulle, W.~E, B.~Engquist, and E.~{Vanden-Eijnden}.
\newblock The heterogeneous multiscale method.
\newblock {\em Acta Numer.}, 21:1--87, 2012.

\bibitem[AH15]{AbdH15}
A.~Abdulle and P.~Henning.
\newblock A reduced basis localized orthogonal decomposition.
\newblock {\em J. Comput. Phys.}, 295:379--401, 2015.

\bibitem[AH17]{AbdH17}
A.~Abdulle and P.~Henning.
\newblock Localized orthogonal decomposition method for the wave equation with
  a continuum of scales.
\newblock {\em Math. Comp.}, 86(304):549--587, 2017.

\bibitem[AHP21]{AltHP21}
R.~Altmann, P.~Henning, and D.~Peterseim.
\newblock Numerical homogenization beyond scale separation.
\newblock {\em Acta Numer.}, 30:1--86, 2021.

\bibitem[BDG20]{BerDG20}
J.~Berner, M.~Dablander, and P.~Grohs.
\newblock Numerically solving parametric families of high-dimensional
  kolmogorov partial differential equations via deep learning.
\newblock In H.~Larochelle, M.~Ranzato, R.~Hadsell, M.~F. Balcan, and H.~Lin,
  editors, {\em Advances in Neural Information Processing Systems}, volume~33,
  pages 16615--16627. Curran Associates, Inc., 2020.

\bibitem[BEKS17]{BezEKS17}
J.~Bezanson, A.~Edelman, S.~Karpinski, and V.~B. Shah.
\newblock Julia: A fresh approach to numerical computing.
\newblock {\em SIAM Rev.}, 59(1):65--98, 2017.

\bibitem[BHKS21]{BhaHKS20}
K.~Bhattacharya, B.~Hosseini, N.~B. Kovachki, and A.~M. Stuart.
\newblock Model reduction and neural networks for parametric {PDEs}.
\newblock {\em SMAI J. Comput. Math.}, 7:121--157, 2021.

\bibitem[BL11]{BabL11}
I.~Babu{\v s}ka and R.~Lipton.
\newblock Optimal local approximation spaces for generalized finite element
  methods with application to multiscale problems.
\newblock {\em Multiscale Model. Simul.}, 9(1):373--406, 2011.

\bibitem[BO00]{BabO00}
I.~Babu{\v{s}}ka and J.~E. Osborn.
\newblock Can a finite element method perform arbitrarily badly?
\newblock {\em Math. Comp.}, 69(230):443--462, 2000.

\bibitem[CE18]{ChaE18}
S.~Chan and A.~H. Elsheikh.
\newblock A machine learning approach for efficient uncertainty quantification
  using multiscale methods.
\newblock {\em J. Comput. Phys.}, 354:493--511, 2018.

\bibitem[CMP20]{CaiMP20}
A.~Caiazzo, R.~Maier, and D.~Peterseim.
\newblock Reconstruction of quasi-local numerical effective models from
  low-resolution measurements.
\newblock {\em J. Sci. Comput.}, 85(1), Article No.~10, 2020.

\bibitem[DG75]{Gio75}
E.~De~Giorgi.
\newblock Sulla convergenza di alcune successioni d'integrali del tipo
  dell'area.
\newblock {\em Rend. Mat. (6)}, 8:277--294, 1975.

\bibitem[EE03]{EE03}
W.~E and B.~Engquist.
\newblock The heterogeneous multiscale methods.
\newblock {\em Commun. Math. Sci.}, 1(1):87--132, 2003.

\bibitem[EE05]{EE05}
W.~E and B.~Engquist.
\newblock The heterogeneous multi-scale method for homogenization problems.
\newblock In {\em Multiscale methods in science and engineering}, volume~44 of
  {\em Lect. Notes Comput. Sci. Eng.}, pages 89--110. Springer, Berlin, 2005.

\bibitem[EG17]{ErnG15}
A.~{Ern} and J.-L. {Guermond}.
\newblock Finite element quasi-interpolation and best approximation.
\newblock {\em ESAIM Math. Model. Numer. Anal.}, 51(4):1367--1385, 2017.

\bibitem[EGH15]{ElfGH15}
D.~Elfverson, V.~Ginting, and P.~Henning.
\newblock On multiscale methods in {P}etrov-{G}alerkin formulation.
\newblock {\em Numer. Math.}, 131(4):643--682, 2015.

\bibitem[EGW11]{EfeGW11}
Y.~R. Efendiev, J.~Galvis, and X.-H. Wu.
\newblock Multiscale finite element methods for high-contrast problems using
  local spectral basis functions.
\newblock {\em J. Comput. Phys.}, 230(4):937--955, 2011.

\bibitem[EH09]{EfeH09}
Y.~R. Efendiev and T.~Y. Hou.
\newblock {\em Multiscale finite element methods}, volume~4 of {\em Surveys and
  Tutorials in the Applied Mathematical Sciences}.
\newblock Springer, New York, 2009.
\newblock Theory and applications.

\bibitem[EHJ17]{WeiHJ17}
W.~E, J.~Han, and A.~Jentzen.
\newblock Deep learning-based numerical methods for high-dimensional parabolic
  partial differential equations and backward stochastic differential
  equations.
\newblock {\em Commun. Math. Stat.}, 5(4):349--380, 2017.

\bibitem[EHJ21]{HanJ20}
W.~E, J.~Han, and A.~Jentzen.
\newblock Algorithms for solving high dimensional {PDEs}: From nonlinear {Monte
  Carlo} to machine learning.
\newblock {\em Nonlinearity}, 35(1):278--310, 2021.

\bibitem[EHMP19]{EngHMP16}
C.~{Engwer}, P.~{Henning}, A.~{M{\aa}lqvist}, and D.~{Peterseim}.
\newblock Efficient implementation of the localized orthogonal decomposition
  method.
\newblock {\em Comput. Methods Appl. Mech. Engrg.}, 350:123--153, 2019.

\bibitem[EY18]{WeiY18}
W.~E and Bing Yu.
\newblock The deep {Ritz} method: a deep learning-based numerical algorithm for
  solving variational problems.
\newblock {\em Commun. Math. Stat.}, 6(1):1--12, 2018.

\bibitem[FP20]{FeiP20}
M.~Feischl and D.~Peterseim.
\newblock Sparse compression of expected solution operators.
\newblock {\em SIAM J. Numer. Anal.}, 58(6):3144--3164, 2020.

\bibitem[GB10]{GloB10}
X.~Glorot and Y.~Bengio.
\newblock Understanding the difficulty of training deep feedforward neural
  networks.
\newblock In {\em Proceedings of the 13th International Conference on
  Artificial Intelligence and Statistics}, pages 249--256. JMLR Workshop and
  Conference Proceedings, 2010.

\bibitem[GHV18]{GalHV18}
D.~Gallistl, P.~Henning, and B.~Verf\"{u}rth.
\newblock Numerical homogenization of {H(curl)}-problems.
\newblock {\em SIAM J. Numer. Anal.}, 56(3):1570--1596, 2018.

\bibitem[GM22]{GeeM21}
S.~Geevers and R.~Maier.
\newblock Fast mass lumped multiscale wave propagation modelling.
\newblock To appear in {\em IMA J. Numer. Anal.}, 2022.

\bibitem[GP15]{GalP15}
D.~Gallistl and D.~Peterseim.
\newblock Stable multiscale {Petrov--Galerkin} finite element method for high
  frequency acoustic scattering.
\newblock {\em Comput. Method. Appl. M.}, 295:1--17, 2015.

\bibitem[GPR{\etalchar{+}}22]{GeiPRSK20}
M.~Geist, P.~Petersen, M.~Raslan, R.~Schneider, and G.~Kutyniok.
\newblock Numerical solution of the parametric diffusion equation by deep
  neural networks.
\newblock To appear in {\em J. Sci. Comput.}, 2022.

\bibitem[GS19]{GhaS19}
F.~Ghavamian and A.~Simone.
\newblock Accelerating multiscale finite element simulations of
  history-dependent materials using a recurrent neural network.
\newblock {\em Comput. Method. Appl. M.}, 357:112594, 2019.

\bibitem[GSW21]{GaoSW21}
H.~Gao, L.~Sun, and J.-X. Wang.
\newblock Phygeonet: physics-informed geometry-adaptive convolutional neural
  networks for solving parameterized steady-state pdes on irregular domain.
\newblock {\em J. Comput. Phys.}, 428:110079, 2021.

\bibitem[HJE18]{HanJW18}
J.~Han, A.~Jentzen, and W.~E.
\newblock Solving high-dimensional partial differential equations using deep
  learning.
\newblock {\em Proceedings of the National Academy of Sciences},
  115(34):8505--8510, 2018.

\bibitem[HJKN20]{HutJKN20}
M.~Hutzenthaler, A.~Jentzen, T.~Kruse, and T.~A. Nguyen.
\newblock A proof that rectified deep neural networks overcome the curse of
  dimensionality in the numerical approximation of semilinear heat equations.
\newblock {\em SN Partial Differ. Equ. Appl.}, 1(2):1--34, 2020.

\bibitem[HKM20]{HelKM20}
F.~Hellman, T.~Keil, and A.~M\aa{}lqvist.
\newblock Numerical upscaling of perturbed diffusion problems.
\newblock {\em SIAM J. Sci. Comput.}, 42(4):A2014--A2036, 2020.

\bibitem[HP13]{HenP13}
P.~Henning and D.~Peterseim.
\newblock Oversampling for the {M}ultiscale {F}inite {E}lement {M}ethod.
\newblock {\em Multiscale Model. Simul.}, 11(4):1149--1175, 2013.

\bibitem[HW97]{Hou1997}
T.~Y. Hou and X.-H. Wu.
\newblock A multiscale finite element method for elliptic problems in composite
  materials and porous media.
\newblock {\em J. Comput. Phys.}, 134(1):169--189, 1997.

\bibitem[HW22]{HenW20}
P.~Henning and J.~W\"arneg\aa{}rd.
\newblock Superconvergence of time invariants for the {Gross-Pitaevskii}
  equation.
\newblock {\em Math. Comp.}, 91(334):509--555, 2022.

\bibitem[HZRS16]{HeZRS16}
K.~He, X.~Zhang, S.~Ren, and J.~Sun.
\newblock Deep residual learning for image recognition.
\newblock In {\em Proceedings of the IEEE conference on computer vision and
  pattern recognition}, pages 770--778, 2016.

\bibitem[Inn18]{Inn18}
M.~Innes.
\newblock Flux: Elegant machine learning with julia.
\newblock {\em J. Open Source Softw.}, 2018.

\bibitem[KB14]{KinB14}
D.~P. Kingma and J.~Ba.
\newblock Adam: A method for stochastic optimization.
\newblock {\em ArXiv Preprint}, 1412.6980, 2014.

\bibitem[KLY21]{KhoLY17}
Y.~Khoo, J.~Lu, and L.~Ying.
\newblock Solving parametric {PDE} problems with artificial neural networks.
\newblock {\em European J. Appl. Math.}, 32(3):421--435, 2021.

\bibitem[KPRS22]{KutPRS19}
G.~Kutyniok, P.~Petersen, M.~Raslan, and R.~Schneider.
\newblock A theoretical analysis of deep neural networks and parametric {PDEs}.
\newblock To appear in {\em Constr. Approx.}, 2022.

\bibitem[MP14]{MalP14}
A.~M{\aa}lqvist and D.~Peterseim.
\newblock Localization of elliptic multiscale problems.
\newblock {\em Math. Comp.}, 83(290):2583--2603, 2014.

\bibitem[MP15]{MalP15}
A.~M\aa{}lqvist and D.~Peterseim.
\newblock Computation of eigenvalues by numerical upscaling.
\newblock {\em Numer. Math.}, 130(2):337--361, 2015.

\bibitem[MP17]{MalP17}
A.~M\aa{}lqvist and D.~Peterseim.
\newblock Generalized finite element methods for quadratic eigenvalue problems.
\newblock {\em ESAIM Math. Model. Numer. Anal.}, 51(1):147--163, 2017.

\bibitem[MP19]{MaiP19}
R.~Maier and D.~Peterseim.
\newblock Explicit computational wave propagation in micro-heterogeneous media.
\newblock {\em BIT Numer. Math.}, 59(2):443--462, 2019.

\bibitem[MP20]{MP20}
A.~M\aa{}lqvist and D.~Peterseim.
\newblock Numerical homogenization by localized orthogonal decomposition.
  volume~5 of {\em SIAM Spotlights}.
\newblock Society for Industrial and Applied Mathematics (SIAM), Philadelphia,
  PA, 2020.
  
\bibitem[MV22]{MV21}
A.~M\aa{}lqvist and B.~Verf{\"u}rth.
\newblock  An offline-online strategy for multiscale problems with random defects.
\newblock To appear in {\em ESAIM Math. Model. Numer. Anal.}, 2022.

\bibitem[MT78]{Mur78}
F.~Murat and L.~Tartar.
\newblock H-convergence.
\newblock {\em S\'{e}minaire d'Analyse Fonctionnelle et Num\'{e}rique de
  l'Universit\'{e} d'Alger}, 1978.

\bibitem[MV22]{MaiV20}
R.~Maier and B.~Verf\"urth.
\newblock Multiscale scattering in nonlinear {K}err-type media.
\newblock To appear in {\em Math. Comp.}, 2022.

\bibitem[Owh17]{Owh17}
H.~Owhadi.
\newblock Multigrid with rough coefficients and multiresolution operator
  decomposition from hierarchical information games.
\newblock {\em SIAM Rev.}, 59(1):99--149, 2017.

\bibitem[OZB14]{OwhZB13}
H.~Owhadi, L.~Zhang, and L.~Berlyand.
\newblock Polyharmonic homogenization, rough polyharmonic splines and sparse
  super-localization.
\newblock {\em ESAIM Math. Model. Numer. Anal.}, 48(2):517--552, 2014.

\bibitem[Pet16]{Pet16}
D.~Peterseim.
\newblock Variational multiscale stabilization and the exponential decay of
  fine-scale correctors.
\newblock In {\em Building Bridges: Connections and Challenges in Modern
  Approaches to Numerical Partial Differential Equations}, volume 114 of {\em
  Lect. Notes Comput. Sci. Eng.}, pages 341--367. Springer, Cham, 2016.

\bibitem[Pet17]{P17}
D.~Peterseim.
\newblock Eliminating the pollution effect in {H}elmholtz problems by local
  subscale correction.
\newblock {\em Math. Comp.}, 86(305):1005--1036, 2017.

\bibitem[PGM{\etalchar{+}}19]{PasG+19}
A.~Paszke, S.~Gross, F.~Massa, A.~Lerer, J.~Bradbury, G.~Chanan, T.~Killeen,
  Z.~Lin, N.~Gimelshein, L.~Antiga, A.~Desmaison, A.~Kopf, E.~Yang, Z.~DeVito,
  M.~Raison, A.~Tejani, S.~Chilamkurthy, B.~Steiner, L.~Fang, J.~Bai, and
  S.~Chintala.
\newblock {PyTorch}: An imperative style, high-performance deep learning
  library.
\newblock In H.~Wallach, H.~Larochelle, A.~Beygelzimer, F.~Alch\'{e}-Buc,
  E.~Fox, and R.~Garnett, editors, {\em Advances in Neural Information
  Processing Systems 32}, pages 8024--8035. Curran Associates, Inc., 2019.

\bibitem[PS12]{PetS12}
D.~Peterseim and S.~A. Sauter.
\newblock Finite {e}lements for {e}lliptic {p}roblems with {h}ighly {v}arying,
  {n}onperiodic {d}iffusion {m}atrix.
\newblock {\em Multiscale Model. Simul.}, 10(3):665--695, 2012.

\bibitem[PS17]{PetS17}
D.~Peterseim and M.~Schedensack.
\newblock Relaxing the {CFL} condition for the wave equation on adaptive
  meshes.
\newblock {\em J. Sci. Comput.}, 72(3):1196--1213, 2017.

\bibitem[PZ21]{PadZ21}
G.~A. Padmanabha and N.~Zabaras.
\newblock A {Bayesian} multiscale deep learning framework for flows in random
  media.
\newblock {\em Found. Data Sci.}, 3(2):251--303, 2021.

\bibitem[RHB19]{RenHB19}
X.~Ren, A.~Hannukainen, and A.~Belahcen.
\newblock Homogenization of multiscale eddy current problem by localized
  orthogonal decomposition method.
\newblock {\em IEEE Transactions on Magnetics}, 55(9):1--4, 2019.

\bibitem[RPK19]{RaiPK19}
M.~Raissi, P.~Perdikaris, and G.~E. Karniadakis.
\newblock Physics-informed neural networks: A deep learning framework for
  solving forward and inverse problems involving nonlinear partial differential
  equations.
\newblock {\em J. Comput. Phys.}, 378:686--707, 2019.

\bibitem[SZ19]{SchZ19}
C.~Schwab and J.~Zech.
\newblock Deep learning in high dimension: Neural network expression rates for generalized polynomial chaos expansions in UQ.
\newblock {\em Analysis and Applications}, 17(01):19--55, 2019.

\bibitem[Spa68]{Spagnolo:1968}
S.~Spagnolo.
\newblock Sulla convergenza di soluzioni di equazioni paraboliche ed
  ellittiche.
\newblock {\em Ann. Sc. Norm. Super. Pisa Cl. Sci.(3)}, 22:571--597, 1968.

\bibitem[SS18]{SirS18}
J.~Sirignano and K.~Spiliopoulos.
\newblock {DGM}: A deep learning algorithm for solving partial differential
  equations.
\newblock {\em J. Comput. Phys.}, 375:1339--1364, 2018.

\bibitem[WCC{\etalchar{+}}20]{WanCCEW20}
Y.~Wang, S.~W. Cheung, E.~T. Chung, Y.~Efendiev, and M.~Wang.
\newblock Deep multiscale model learning.
\newblock {\em J. Comput. Phys.}, 406:109071, 2020.

\end{thebibliography}
\end{document}